\documentclass[reqno,11pt]{amsart}

\usepackage{color}
\usepackage{graphicx}
\usepackage[hidelinks]{hyperref}
\usepackage{enumerate}
\usepackage{latexsym,array}
\usepackage{amsfonts}
\usepackage{shadow}
\usepackage{amsbsy}
\usepackage{dsfont}

\newcommand{\id}{\mathcal{I}}
\newcommand{\ext}{\mathrm{ext}}
\newcommand{\bimon}{\mathcal{N}}
\newcommand{\lin}{\mathrm{span}}

\newtheorem{Pa}{Paper}[section]
\newtheorem{Tm}[Pa]{{\bf Theorem}}
\newtheorem{La}[Pa]{{\bf Lemma}}
\newtheorem{Dn}[Pa]{{\bf Definition}}
\newtheorem{Cy}[Pa]{{\bf Corollary}}

\newtheorem{Pn}[Pa]{{\bf Proposition}}

\newtheorem{Ex}[Pa]{{\bf Example}}

\theoremstyle{definition}
\newtheorem{Rk}[Pa]{{\bf Remark}}

\def\C{\mathbb C}
\def\hh{\mathbb{H}}

\newcommand{\pp}{\partial}
\newcommand{\rr}{\mathbb{R}}
\renewcommand{\S}{\mathbb{S}}

\newcommand{\boundOP}{\mathcal{B}}

\newcommand{\SM}{\mathcal{SM}}

\renewcommand{\Re}{\mathrm{Re}}

\newcommand{\prodL}{{\star_\ell}}
\newcommand{\prodR}{{\star_r}}

\newcommand{\dist}{\mathrm{dist}}

\title[
On power series expansions of the S-resolvent operator and the Taylor formula]
{On power series expansions of the S-resolvent operator and the Taylor formula} \oddsidemargin
0.2in \evensidemargin 0.2in \topmargin -0.5in \textwidth
15.5truecm \textheight 23truecm
\def\H{\mathbb H}

\def\R{\mathbb R}
\def\N{\mathbb N}
\def\C{\mathbb C}

\def\Z{\mathbb Z}

\def\(s){\mathscr S(\R\times\R)}

\author[F. Colombo]{Fabrizio Colombo}
\address{(FC)
Politecnico di Milano\\Dipartimento di Matematica\\Via E. Bonardi, 9\\20133
Milano, Italy}
\email{fabrizio.colombo@polimi.it}

\author[J. Gantner]{Jonathan Gantner}
\address{(JG)
Politecnico di Milano\\Dipartimento di Matematica\\Via E. Bonardi, 9\\20133
Milano, Italy
} \email{jonathan.gantner@gmx.at}


\begin{document}
\maketitle

\parindent 0cm
\begin{abstract}
 The  $S$-functional calculus is based on the theory of slice hyperholomorphic functions and it
 defines functions of $n$-tuples of not necessarily commuting operators or of quaternionic operators.
 This calculus relays on the notion of $S$-spectrum and of
  $S$-resolvent operator. Since most of the properties that hold for  the Riesz-Dunford functional calculus extend to the S-functional calculus
  it can be considered its  non commutative version.
In this paper we show that the Taylor formula of  the Riesz-Dunford functional calculus can be generalized  to the S-functional calculus, the proof is not a trivial extension of the classical case because there are several obstructions
due to the non commutativity of the setting in which we work
that have to be overcome.
To prove the Taylor formula we need to introduce a new series expansion of the $S$-resolvent operators associated to the sum of
two $n$-tuples of operators. This result is a crucial step in the proof of our main results,
 but it is also of independent interest because it gives a new series expansion for the $S$-resolvent operators.
This paper is devoted to researchers working in operators theory and hypercomplex analysis.
\end{abstract}

\parindent 0cm



\section{Introduction}

The theory of slice hyperholomorphic functions, introduced in the last decade, turned out to be an important tool in operator theory.
We call slice hyperholomorphic functions two classes of hyperholomorphic functions: the slice monogenic and the slice regular ones.
Slice monogenic functions are defined on subsets of the Euclidean space $\mathbb{R}^{n+1}$ and have values in the real Clifford algebra $\mathbb{R}_n$, while slice regular functions  are defined on subsets of the quaternions and are quaternion-valued.
 These two sets of functions are in the kernel of a suitable
 Cauchy-Riemann operator and so they constitute the class of hyperholomorphic functions for which the S-functional calculus is defined.
 This class of functions plays the analogue role of the holomorphic functions for the Riesz-Dunford functional calculus.

Precisely, when the S-functional calculus is based on slice monogenic functions, it defines functions of $n$-tuples of not necessarily commuting operators, see
\cite{CONVACS, acgs, cs, functionalcss},
its commutative version is in \cite{SCFUNCTIONAL}, while when
we consider slice regular functions the S-functional calculus defines functions of quaternionic operators (the quaternionic functional calculus),
 see \cite{acgs, JGA,  formulation}.
 To define the S-functional calculus a new notion of spectrum for $n$-tuples of operators and for quaternionic operators has been introduced.
 It is called S-spectrum and it naturally arise from the Cauchy formula of slice hyperholomorphic functions.
 For a global picture of the function theory and of the S-functional calculus see the book \cite{css_book}.

There exists a continuous functional calculus for quaternionic bounded linear normal operators in a Hilbert space which is
based on the notion of $S$-spectrum and it has been considered in \cite{GMP}.
 On the S-spectrum is also based a spectral theorem for quaternionic unitary operators see \cite{acks2} whose proof uses a theorem of Herglotz in  \cite{acks1}. For the spectral theorem of compact quaternionic normal operators see \cite{spectcomp}.
 For the spectral theorem based on the $S$-spectrum for quaternionic bounded or unbounded normal operators see the paper \cite{Spectralnorm}.

In this paper we will formulate our results for the case of $n$-tuples of non commuting operators,
 but what we prove here holds also for the quaternionic functional calculus.

Let us recall the Taylor formula for  the Riesz-Dunford functional calculus, see \cite{ds} p.590 for more details,
 and then we show its generalization to our setting.

{\it Let $A$ and $N$ be bounded commuting operators on a complex Banach space $X$.
Let $f$ be an analytic function on a domain $D\subset \mathbb{C}$ that includes the spectrum $\sigma(A)$ of $A$
and every point within a distance of $\sigma(A)$ not greater than some positive number $\varepsilon$.
Suppose that the spectrum $\sigma(N)$ of $N$ lies within the open circle of radius $\varepsilon$
about the origin.

Then $f$ is analytic on a neighborhood of $\sigma(A+N)$ and the series
$$
f(A+N)=\sum_{m= 0}^{\infty}\frac{f^{(m)}(A)}{m!}N^m,
$$
 converges in the uniform operator topology.}

One of the main points in the proof of this theorem is the series expansion of the resolvent operator
$R(\lambda, A+N):=(\lambda I-A-N)^{-1}$ in terms of the resolvent of $A$ and in powers of $N$. Precisely the series
\begin{equation}\label{VAN}
V_\lambda(A,N):=\sum_{m= 0}^{\infty}R(\lambda, A)^{m+1}N^m
\end{equation}
converges uniformly for $\lambda$ in a set whose minimum distance from $\sigma(A)$ is greater than $\varepsilon$, and since $A$ and $N$ commute it is easy to verify that
$$
V_\lambda(A,N)(\lambda I-A-N)=(\lambda I-A-N)V_\lambda(A,N) = \id,
$$
so $V_\lambda(A,N)$ is the required power series expansion for the resolvent operator of $A+N$.
\\
\\
In the case of $n$-tuples of not necessarily commuting operators we will see that the analogue of (\ref{VAN}) is not so easy to find, but despite this difficulty we are still able to prove the Taylor formula.

Let $T_j:V\to V$, $j=0,1,...,n$ be bounded $\mathbb{R}$-linear operators defined on a real Banach space~$V$. Associated with this $(n+1)$-tuple of operators
we define the paravector operator $T$ as  $T:=T_0+e_1T_1+\ldots+T_n$ where $e_j$, $j=1,\ldots,n$ are the units of the real Clifford algebra $\mathbb{R}_n$.
 In this paper it will always consider bounded paravector operators.

The $S$-spectrum of $T$ is
$$
\sigma_S(T):=\{ s\in \mathbb{R}^{n+1}\ \ :\ \ T^2-2 \Re(s)T+|s|^2\mathcal{I}\ \ \
{\it is\ not\  invertible}\}.
$$
The $S$-resolvent set is defied as $\rho_S(T)= \mathbb{R}^{n+1}\setminus  \sigma_S(T)$,
the left $S$-resolvent operator is
$$
S_L^{-1}(s,T):=-(T^2-2\Re(s) T+|s|^2\mathcal{I})^{-1}(T-\overline{s}\mathcal{I}),\ \ \ \ s\in \rho_S(T)
$$
and the right $S$-resolvent operator is
$$
S_R^{-1}(s,T):=-(T-\overline{s}\mathcal{I})(T^2-2\Re(s) T+|s|^2\mathcal{I})^{-1},\ \ \ \ s\in \rho_S(T).
$$
In this case the analogue of the power series (\ref{VAN}) requires a new series expansion for the  S-resolvent operators
which involves,  for $n\geq 0$, the operator
\[
S_L^{-n}(s,T):= (T^2 - 2s_0T + |s|^2\id)^{-n}(\overline{s}\:\!\id-T)^{\prodL n},\ \ \ \ s\in \rho_S(T)
\]
where the Newton binomial $(\overline{s}\:\!\id-T)^{\prodL n}$ is computed allowing $\overline{s}$ and $T$ to commute in a sense
explained in the sequel.
 When we assume that the bounded paravector operators $T$ and  $N$ commute  and
  are such that $\sigma_S(N)\subset B_{\varepsilon}(0)$. Then, for any $s\in\rho_S(T)$ with $\dist(s,\sigma_S(T)) > \varepsilon$,
 we obtain the following expansions for the S-resolvent operators:
 \[
 S_L^{-1}(s,T+N) = \sum_{n=0}^\infty N^n\, S_L^{-(n+1)}(s,T)
 \qquad\text{and}\qquad S_R^{-1}(s,T+N) = \sum_{n=0}^{\infty}S_R^{-(n+1)}(s,T)\, N^n.
 \]
 Moreover, these series converge uniformly on any set $C$ with ${\dist(C,\sigma_S(T))>\varepsilon}$.

A crucial fact to obtain these series expansions is to prove the invertibility of the term
$$(T+N)^2 - 2s_0(T+N) + |s|^2\id$$ under the required  conditions on the operators $T$ and $N$.

Precisely, in Theorem \ref{SpectrumStable}, it is proved that
if  $T$ and  $N$ are bounded paravector operators  that commute and are such that $\sigma_S(N)$ is contained in the open ball $B_{\varepsilon}(0)$, then, when  $\dist(s,\sigma_S(T)) > \varepsilon$, we have that $s\in\rho_S(T+N)$ and
\begin{equation*}
\left((T+N)^2 - 2s_0(T+N) + |s|^2\id\right)^{-1} = \sum_{n=0}^\infty\left(\sum_{k=0}^nS_L^{-(k+1)}(s,T)\prodL S_L^{-(n-k+1)}(\overline{s},T)\right)N^n,
\end{equation*}
where the series converges in the operator norm.

With the above results we are the able to prove the Taylor formula for the $S$-functional calculus:
under the conditions that ${T, N\in\boundOP^{0,1}(V_n)}$ are  such that $T$ and $N$ commute and  $\sigma_S(N)\subset B_{\varepsilon}(0)$,
 given a  slice hyperholomorphic function $f$ defined on some suitable domain $U$ that contains $\sigma_S(T)$
  and every point $s\in\R^{n+1}$ with $\dist( s,\sigma_S(T))\leq \varepsilon$, we have
\[
f(T+N) =  \sum_{n=0}^{\infty} N^n\frac{1}{n!} \left(\partial_{s}^nf\right)(T).
 \]
An analogous formula holds for right slice hyperholomorphic functions.
The Taylor formula can be applied, for example, in the case that  $N$ is a Riesz projector of the $S$-functional calculus
associated with $T$, because in this case $TN=NT$.

The  $S$-resolvent operators defined above are also important in the realization of Schur functions in the slice hyperholomorphic setting,
see \cite{aacks, ACLS, acs1, acs2, acs3}.
For  Schur analysis in the classical case see \cite{MR2002b:47144, adrs}.
\\
\\
The plan of the paper is as follows.
In Section 2 we recall the main results on slice monogenic functions that are useful in the sequel.
Here we also recall the notion of holomorphic maps of a
paravector variable which allows us to consider axially symmetric
domains that do not necessarily intersect the real line.
This is necessary in order to study our problem without
restrictions on the set $C$ such that  $\dist(C,\sigma_S(T)) > \varepsilon$, (see the introduction).
In Section~3 we recall the main results on the $S$-functional calculus that we will use in the sequel.
In Section~4 we give an alternative proof of the invertibility of the operator $T^2-2\Re(s)T+|s|^2\mathcal{I}$
for $\|T\|< |s|$ using the product of the two power series $\sum_{n= 0}^\infty T^ns^{-1-n}$
and $\sum_{n=0}^\infty (\overline{s}^{-1}T)^n\overline{s}^{-1}$ which both converge for $\|T\|< |s|$.
A more direct proof first appeared in \cite{GMP} in the quaternionic setting.
In Section 5 we consider slice hyperholomorphic powers of the $S$-resolvent operators that we will
use in Section 6 in order to prove the Taylor formula for the $S$-functional calculus in Theorem \ref{THTAYLOR}.

This paper is devoted to researchers working in operators theory and in hypercomplex analysis, for this reason we recall the main
results in both research areas in order to make it accessible to both communities.

\section{Preliminary results in function theory}

We advise the reader that in the following the symbol $n$ is a natural number that indicates the number of imaginary units of the real Clifford algebra $\mathbb{R}_n$,  the dimension of the Euclidean space $\mathbb{R}^{n+1}$, but it will be used also  as an index. The context makes clear the meaning of the symbol.
In this section we recall the fundamentals of the theory of slice hyperholomorphic functions,  see \cite{css_book}.
\\
Let $\rr_n$ be the real Clifford algebra over $n$ imaginary units $e_1,\ldots ,e_n$
satisfying the relations $e_ie_j+e_je_i=0$ for  $i\not= j$ and  $e_i^2=-1$ or, equivalently,
\begin{equation}\label{DefRel} e_ie_j + e_je_i = - 2\delta_{ij}.\end{equation}
 An element in the Clifford algebra will be denoted by $\sum_A e_Ax_A$ where
$A=\{ i_1\ldots i_r\}\in \mathcal{P}\{1,2,\ldots, n\},\ \  i_1<\ldots <i_r$
 is a multi-index
and $e_A=e_{i_1} e_{i_2}\ldots e_{i_r}$, $e_\emptyset =1$.
An element $(x_0,x_1,\ldots,x_n)\in \mathbb{R}^{n+1}$  will be identified with the element
$
 x=x_0+\underline{x}=x_0+ \sum_{j=1}^nx_je_j\in\mathbb{R}_n,
$
which is called a paravector. The real part $x_0$ of $x$ will also be denoted by $\Re(x)$.
The norm of $x\in\mathbb{R}^{n+1}$ is defined by $|x|^2=x_0^2+x_1^2+\ldots +x_n^2$
 and the conjugate of $x$ is defined as
$
\bar x=x_0-\underline x=x_0- \sum_{j=1}^nx_je_j.
$
Let
$$
\S:=\{ \underline{x}=e_1x_1+\ldots +e_nx_n\ | \  x_1^2+\ldots +x_n^2=1\}.
$$
For $I\in\S$, we obviously have $I^2=-1$.
Given an element $x=x_0+\underline{x}\in\rr^{n+1}$, we set
\[
I_x:=\begin{cases}\underline{x}/|\underline{x}| & \text{if }\underline{x}\neq0\\
\text{any element of $\S$} & \text{if } \underline{x} = 0.\end{cases}
\]
Then $x = u + I_xv$ with $u = x_0$ and $v = |\underline{x}|$.

For any element $x=u + I_xv\in\rr^{n+1}$, the set
$$
[x]:=\{y\in\rr^{n+1}\ :\ y=u+I v, \ I\in \mathbb{S}\}
$$

is an $(n-1)$-dimensional sphere in $\mathbb{R}^{n+1}$.
The vector space $\mathbb{R}+I\mathbb{R}$ passing through $1$ and
$I\in \mathbb{S}$ will be denoted by $\mathbb{C}_I$ and
an element belonging to $\mathbb{C}_I$ will be indicated by $u+Iv$ with  $u$, $v\in \mathbb{R}$.

Since we identify the set of paravectors with the space $\R^{n+1}$, if $U\subseteq\mathbb{R}^{n+1}$ is an open set,
a function $f:\ U\subseteq \mathbb{R}^{n+1}\to\mathbb{R}_n$ can be interpreted as
a function of a paravector $x$.

\begin{Dn}[Slice hyperholomorphic functions]
\label{defsmon}
Let $U\subseteq\mathbb{R}^{n+1}$ be an open set and let
$f: U\to\mathbb{R}_n$ be a real differentiable function. Let
$I\in\mathbb{S}$ and let $f_I$ be the restriction of $f$ to the
complex plane $\mathbb{C}_I$.
\\
The function  $f$ is said to be left slice  hyperholomorphic (or slice monogenic) if, for every
$I\in\mathbb{S}$, it satisfies
$$
\frac{1}{2}\left(\frac{\partial }{\partial u}f_I(u+Iv)+I\frac{\partial
}{\partial v}f_I(u+Iv)\right)=0
$$
 on $U\cap \mathbb{C}_I$. We denote the set of left slice hyperholomorphic functions on the open set $U$  by $\mathcal{SM}^L(U)$.
\\
The function $f$ is said to be right slice hyperholomorphic (or right slice monogenic) if,
for every
$I\in\mathbb{S}$, it satisfies
$$
\frac{1}{2}\left(\frac{\partial }{\partial u}f_I(u+Iv)+\frac{\partial
}{\partial v}f_I(u+Iv)I\right)=0
$$
on $U\cap \mathbb{C}_I$.
We denote the set of right slice hyperholomorphic functions on the open set $U$  by $\mathcal{SM}^R(U)$.
\end{Dn}
Any power series of the form $\sum_{n=0}^\infty x^n a_n$ with $a_\ell\in\mathbb{R}_n$, for $\ell\in \mathbb{N}$, is left slice hyperholomorphic and any power series of the form $\sum_{n=0}^\infty b_nx^n$ with $b_\ell\in\mathbb{R}_n$, for $\ell\in \mathbb{N}$, is right slice hyperholomorphic. Conversely, at any real point, any left or right slice hyperholomorphic function allows a power series expansion of the respective form.
\begin{Dn}\label{SliceDeriv}
Let $U\subset\R^{n+1}$ be open. The slice derivative of a function $f\in\SM^L(U)$ is the function defined by
\[ \partial_s f (u + Iv) := \frac12\left(\frac{\partial}{\partial u} f_I(u + Iv) - I\frac{\partial}{\partial v}f_I(u + Iv)\right)\]
for $u + I v\in U$.

The slice derivative of a function $f\in\SM^R(U)$ is the function defined by
\[ \partial_s f (u + Iv) := \frac12\left(\frac{\partial}{\partial u} f_I(u + Iv) - \frac{\partial}{\partial v}f_I(u + Iv)I\right)\]
for $u + I v \in U$.
\end{Dn}
\begin{La}
Let $U\subset\R^{n+1}$ be open. If $f\in\SM^L(U)$, then $\partial_s f\in\SM^L(U)$. If $f\in\SM^R(U)$, then $\partial_s f\in\SM^R(U)$.
\end{La}
\begin{Rk}\label{SliceDerivS0}
Note that the slice derivative coincides with the partial derivative with respect to the real part since
\begin{align*}\partial_s f (u+ Iv) &= \frac12\left(\frac{\partial}{\partial u} f_I(u + Iv) - I\frac{\partial}{\partial v}f_I(u + Iv)\right) \\
&= \frac12\left(\frac{\partial}{\partial u} f_I(u + Iv) + \frac{\partial}{\partial u}f_I(u + Iv)\right) = \frac{\partial}{\partial u}f_I(u + Iv)\end{align*}
for any left slice hyperholomorphic function. An analogous computation shows the respective identity for right slice hyperholomorphic functions.
\end{Rk}
We denote by $B_r(y)$ the open ball in $\R^{n+1}$ of radius $r>0$ and centered at $y\in \mathbb{R}^{n+1}$.
\begin{La}
Let $y = y_0 \in\R$. A function $f: B_r(y)\subset\R^{n+1}\to\R_n$ is left slice hyperholomorphic if and only if
\[f (x) = \sum_{n=0}^{\infty} x^n \frac{1}{n!}\partial_s^nf(y). \]
A function $f: B_r(y)\subset\R^{n+1}\to\R_n$ is right slice hyperholomorphic if and only if
\[f (x) = \sum_{n=0}^{\infty} \frac{1}{n!}\partial_s^nf(y) \,x^n . \]
\end{La}

\begin{La}[Splitting Lemma]\label{SplitLem}
Let $U\subset\R^{n+1}$ be open and let $f:U\to\R_n$ be a left slice hyperholomorphic function. For every $I=I_1\in\S$, let $I_2,\ldots,I_n$ be a completion to a basis of $\R_n$ that satisfies the relation \eqref{DefRel}, that is $I_rI_s + I_s I_r = - 2\delta_{r,s}$. Then there exist $2^{n-1}$ holomorphic functions $F_A: U\cap\C_I\to\C_I$ such that for every $z\in U\cap\C_I$
\[f_I (z) = \sum_{|A|=0}^{n-1}F_A(z)I_A\]
where $I_A = I_{i_1}\cdots I_{i_s}$ for any nonempty subset $A=\{i_1<\ldots<i_s\}$  of $\{2,\ldots,n\}$ and $I_{\emptyset} = 1$.

Similarly, if $g:U\to\H$ is right slice hyperholomorphic, then there exist $2^{n-1}$ holomorphic functions $G_A:U\cap\C_I\to\C_I$ such that for every $z\in U\cap\C_I$
\[g_I(z) = \sum_{|A|=0}^{n-1}I_AG_A(z).\]
\end{La}

Slice hyperholomorphic functions possess good properties when they are defined on suitable domains which are
introduced in the following definition.

\begin{Dn}[Axially symmetric slice domain]\label{axsymm}
Let $U$ be a domain in $\rr^{n+1}$.
We say that $U$ is a
\textnormal{slice domain} (s-domain for short) if $U \cap \mathbb{R}$ is nonempty and if $U\cap \mathbb{C}_I$ is a domain in $\mathbb{C}_I$ for all $I \in \mathbb{S}$.
We say that $U$ is
\textnormal{axially symmetric} if, for all $x \in U$, the
$(n-1)$-sphere $[x]$ is contained in $U$.
\end{Dn}

\begin{Tm}[The Cauchy formula with slice hyperholomorphic kernel]
\label{Cauchygenerale}
Let $U\subset\mathbb{R}^{n+1}$ be an axially symmetric s-domain.
Suppose that $\partial (U\cap \mathbb{C}_I)$ is a finite union of
continuously differentiable Jordan curves  for every $I\in\mathbb{S}$ and set  $ds_I=-ds I$ for $I\in \mathbb{S}$.
\begin{itemize}
\item
If $f$ is
a (left) slice hyperholomorphic function on a set that contains $\overline{U}$ then
\begin{equation}\label{Cauchyleft}
 f(x)=\frac{1}{2 \pi}\int_{\partial (U\cap \mathbb{C}_I)} S_L^{-1}(s,x)\,ds_I\, f(s)
\end{equation}
where
\begin{equation}\label{SL1}
S_L^{-1}(s,x):=-(x^2 -2 \Re(s)x+|s|^2)^{-1}(x-\overline s).
\end{equation}
\item
If $f$ is a right slice hyperholomorphic function on a set that contains $\overline{U}$,
then
\begin{equation}\label{Cauchyright}
 f(x)=\frac{1}{2 \pi}\int_{\partial (U\cap \mathbb{C}_I)}  f(s)\,ds_I\, S_R^{-1}(s,x)
 \end{equation}
 where
 \begin{equation}\label{SR1}
S_R^{-1}(s,x):= -(x-\bar s)(x^2-2{\rm Re}(s)x+|s|^2)^{-1}.
\end{equation}
\end{itemize}
The
integrals  depend neither on $U$ nor on the imaginary unit
$I\in\mathbb{S}$.
\end{Tm}
The above Cauchy formulas are the  starting point to define the S-functional calculus.
A crucial fact about slice hyperholomorphic functions is the following Representation Formula (also called  Structure Formula).

\begin{Tm}[Representation Formula]\label{formulaMON} Let $U\subset\H$ be an axially symmetric s-domain and let $I\in\S$.
\begin{itemize}
\item
If $f$ is a left slice hyperholomorphic function on $U$, then the following identity holds true for all $x=u+Jv \in U $:
\begin{equation}\label{distribution}
f(u+Jv) =\frac{1}{2}\Big[   f(u+Iv)+f(u-Iv)\Big] +J\frac{1}{2}\Big[ I[f(u-Iv)-f(u+Iv)]\Big].
\end{equation}
Moreover, for all $u, v \in \mathbb{R}$ with $u+Iv\in U$, the quantities
\begin{equation}\label{cappa}
\alpha(u,v):=\frac{1}{2}\Big[   f(u+Iv)+f(u-Iv)\Big] \ \ {\sl and}
\ \  \beta(u,v):=\frac{1}{2}\Big[ I[f(u-Iv)-f(u+Iv)]\Big]
\end{equation}
are independent of the imaginary unit $I\in\S$.
\item
If $f$ is a right slice hyperholomorphic function on $U$, then the following identity holds true for all $x=u+Jv \in U $:
\begin{equation}\label{distributionright}
f(u+Jv) =\frac{1}{2}\Big[   f(u+Iv)+f(u-Iv)\Big] +\frac{1}{2}\Big[ [f(u-Iv)-f(u+Iv)]I\Big]J.
\end{equation}
Moreover, for all $u, v \in \mathbb{R}$ with $u+Iv\in U$, the quantities
\begin{equation}\label{capparight}
\alpha(u,v):=\frac{1}{2}\Big[   f(u+Iv)+f(u-Iv)\Big] \ \ {\sl and}
\ \  \beta(u,v):=\frac{1}{2}\Big[ [f(u-Iv)-f(u+Iv)]I\Big]
\end{equation}
are independent of the imaginary unit $I\in\S$.
\end{itemize}
\end{Tm}

\begin{La}[Extension Lemma]
Let $I\in\S$ and let $U$ be a domain in $\C_I$ that is symmetric with respect to the real axis. Then the set
\[[U] = \{ u + Jv: u+ Iv \in U, J\in\S  \}\]
is an axially symmetric s-domain.

Any function $f: U\to\hh$ that satisfies ${\frac{\partial}{\partial u} f (u+Iv) + I\frac{\partial}{\partial v} f(u+Iv) = 0}$ has a unique left slice hyperholomorphic extension $\ext_L(f)$ to $[U]$, which is determined by
\[\ext_L (f )( u + Jv) = \frac{1}{2}\Big[   f(u+Iv)+f(u-Iv)\Big] +J\frac{1}{2}\Big[ I[f(u-Iv)-f(u+Iv)]\Big] .\]
Any function $g: U\to\hh$ that satisfies $\frac{\partial}{\partial u} f (u+Iv) + \frac{\partial}{\partial v} f(u+Iv) I= 0$ has a unique right slice hyperholomorphic extension $\ext_R(f)$ to $[U]$, which is determined by
\[\ext_R(f)(u+Jv) =\frac{1}{2}\Big[   f(u+Iv)+f(u-Iv)\Big]+\frac{1}{2}\Big[ [f(u-Iv)-f(u+Iv)]I\Big]J.\]
\end{La}
The product of two slice hyperholomorphic functions is in general not slice hyperholomorphic. We introduce now a regularized product that preserves the slice hyperholomorphicity.

\begin{Dn}
Let $f$ and $g$ be left slice hyperholomorphic functions on an axially symmetric slice domain. Let $I\in\S$ and let $f_I(z) = \sum_{|A|=0}^{n} F_A(z)I_A$ and $g_I(z) = \sum_{|B|=0}^{n} G _B(z)I_B$ as in Lemma~\ref{SplitLem}. We set
\[f_I\prodL g_I (z) :=  \sum_{|A|\ \mathrm{even}, B} (-1)^{\frac{|A|}{2}}F_A(z)G_B(z)I_AI_B + \sum_{|A|\ \mathrm{odd}, B}(-1)^{\frac{|A|+1}{2}}F_A(z)\overline{G_B(\overline{z})}I_AI_B
\]
for $z\in U\cap\C_I$ and we define the left slice hyperholomorphic product of $f$ and $g$ as
\[f\prodL g := \ext_L(f_I\prodL g_I).\]

For two right slice hyperholomorphic functions $f$ and $g$, we write $f_I(z) = \sum_{|A|=0}^{n} I_AF_A(z)$ and $g_I(z) = \sum_{|B|=0}^{n} I_BG _B(z)$ as in Lemma~\ref{SplitLem} and set
\[f_I\prodR g_I (z) :=  \sum_{A, |B|\ \mathrm{even}}  (-1)^{\frac{|B|}{2}}I_AI_BF_A(z)G_B(z) + \sum_{A, |B|\ \mathrm{odd}}(-1)^{\frac{|B|+1}{2}}I_AI_B \overline{F_A(\overline{z})}G_B(z)
\]
for $z\in U\cap\C_I$ and we define their right slice hyperholomorphic product as
\[f\prodR g := \ext_R(f_I\prodR g_I).\]
\end{Dn}
\begin{Rk}\label{RemarkSimpleProduct}Let $f,g\in\SM^L(U)$ and let $I\in\S$. If $f_I(U\cap\C_I)\subset\C_I$, then $F_A = 0$ if $|A|>0$. Hence, $(f_I\prodL g_I)(x) = f_I(x) g_I(x)$ for any $x\in U\cap\C_I$. However, note that this does not necessarily imply $(f\prodL g)(x) = f(x) g(x)$ for all $x\in U$ because, in general, $\ext_L(f_Ig_I)\neq f g$.

Similarly, if $f,g\in\SM^R(U)$ and $g_I(U\cap\C_I)\subset\C_I$, then $(f\prodR g)(x) = f(x)g(x)$ for any $x\in U\cap\C_I$.
\end{Rk}
\begin{La} \label{ProdProp}Let $U\subset\R^{n+1}$ be an axially symmetric slice domain.
\begin{enumerate}[(i)]
\item If $f,g\in\SM^L(U)$, then $f\prodL g$ is left slice hyperholomorphic and does not depend on the imaginary unit $I$ used in its definition. If $f,g\in\SM^R(U)$, then $f\prodR g$ is right slice hyperholomorphic and does not depend on the imaginary unit $I$ used in its definition.
\item The left and the right slice hyperholomorphic product are associative and distributive over the pointwise addition.
\item \label{ProdSeries}For two left slice hyperholomorphic power series $f(x) = \sum_{n=0}^{\infty} x^na_n$ and $g(x) = \sum_{n=0}^{\infty}x^nb_n$, we have
\begin{equation}\label{prodLSeries}f\prodL g(x) = \sum_{n=0}^{\infty}x^n\left(\sum_{k=0}^na_kb_{n-k}\right).\end{equation}
For two right slice hyperholomorphic power series $f(x) = \sum_{n=0}^{\infty}a_n x^n$ and $g(x) = \sum_{n=0}^{\infty}b_nx^n$, we have
\begin{equation}\label{prodRSeries}f\prodR g(x) = \sum_{n=0}^{\infty}\left(\sum_{k=0}^na_kb_{n-k}\right)x^n.\end{equation}

\end{enumerate}
\end{La}

\begin{Dn}
Let $U\subset\rr^{n+1}$ be an axially symmetric s-domain and let $I\in\S$. For ${f\in\SM^L(U)}$, we write $f_I(z) = \sum_{|A| = 0}^{n-1}F_A(z)I_A$ as in Lemma~\ref{SplitLem} and define
\[ f_I^c(z) := \sum_{|A|\equiv 0} \overline{F_A(\overline{z})}I_A- \sum_{|A| \equiv1}F_A(z)I_A - \sum_{|A|\equiv 2} \overline{F_A(\overline{z})}I_A+ \sum_{|A| \equiv3}F_A(z)I_A, \]
where the equivalence $\equiv$ is intended as $\equiv (\mathrm{mod} 4)$, i.e., the congruence modulo $4$. We define the left slice hyperholomorphic conjugate of $f$ as
\[f^c := \ext_L (f_I^c)\]
and we define the left slice hyperholomorphic inverse of $f$ as
\[f^{-\prodL} := (f^c\prodL f)^{-1} f^c.\]

 For $g\in\SM^R(U)$, we write $g_I(z) = \sum_{|A| = 0}^{n-1}I_AG_A(z)$ as in Lemma~\ref{SplitLem} and define
\[
g_I^c(z) := \sum_{|A|\equiv 0} I_A\overline{G_A(\overline{z})}- \sum_{|A| \equiv1}I_AG_A(z) - \sum_{|A|\equiv 2}I_A \overline{G_A(\overline{z})}+ \sum_{|A| \equiv3}I_AG_A(z).
 \]
 We define the right slice hyperholomorphic conjugate of $g$ as
\[g^c := \ext_R(g_I^c)\]
and we define the right slice hyperholomorphic inverse of $g$ as
\[g^{-\prodL} :=  g^c(g^c\prodR g)^{-1}.\]

\end{Dn}

\begin{La}
Let $U$ be an axially symmetric s-domain.
\begin{enumerate}[(i)]
\item If $f\in\SM^L(U)$, then $f^c$ is left slice hyperholomorphic on $U$ and does not depend on the imaginary unit used in its definition.  If $g\in\SM^R(U)$, then $g^c$ is right slice hyperholomorphic on $U$ and does not depend on the imaginary unit used in its definition.
\item If $f\in\SM^L(U)$ then $f^{-\prodL}$ is left slice hyperholomorphic on the set
$ U\setminus[\mathcal{Z}(f)],$
where $\mathcal{Z}(f) = \{ x\in U :  f(x) = 0\}$ and $[\mathcal{Z}(f)] = \bigcup_{x\in\mathcal{Z}(f)}[x]$. Moreover,
\[f^{-\prodL}\prodL f = f\prodL f^{-\prodL} = 1.\]
If $g\in\SM^R(U)$ then $g^{-\prodR}$ is left slice hyperholomorphic on the set
$ U\setminus[\mathcal{Z}(g)]$ and
\[g^{-\prodR}\prodR g = g\prodR g^{-\prodR} = 1.\]
\end{enumerate}
\end{La}
\begin{Dn}
Let $U\subset \R^{n+1}$ be open. We define
\[ \bimon(U) := \left\{f \in\SM^L(U): f(U\cap\C_I)\subset\C_I\quad\forall I\in\S\right\}.\]
\end{Dn}
The class $\bimon(U)$ plays a privileged role within the set of left slice hyperholomorphic functions. It contains all power series with real coefficients and any function that belongs to it is both left and right slice hyperholomorphic.

\pagebreak[2]
\begin{La}\label{BimonProp}Let $U$ be an axially symmetric slice domain.
\begin{enumerate}[(i)]
\item \label{BimonProd} If $f\in\bimon(U)$, $g\in\SM^L(U)$ and $h\in\SM^R(U)$, then
\[f(x)g(x) = (f\prodL g)(x) = (g\prodL f) (x)
\qquad\text{and}\qquad h(x)f(x) = (h\prodR f) (x) = (f\prodR h)(x).\]
\item If $f\in \bimon(U)$, then $f^c = f$ and $f^{-\prodL} = f^{-\prodR} = f^{-1}$.
\item If $f\in \SM^L(U)$, then $f\prodL f^c\in \bimon(U)$.  If $f\in \SM^R(U)$, then $f\prodR f^c\in \bimon(U)$.
\end{enumerate}
\end{La}
\begin{Ex}\label{ExSliceInv}
Let $s\in\R^{n+1}$. We consider the function $f(x) =  s - x$, which is both
 left and right slice hyperholomorphic. Then $f^c(x) = \overline{s - \overline{x}}= \overline{s} - x$,  and in turn $(f\prodL f^c) (x) = x^2 - 2(s+\overline{s})x + s\overline{s} = x^2 - 2s_0x + |s|^2$ by (\ref{ProdSeries}) in Lemma~\ref{ProdProp}. This function belongs to $\bimon(\R^{n+1})$ because its coefficients are real, and hence $x\in\C_I$ implies $f\prodL f^c(x)\in \C_I$. Moreover, the left slice hyperholomorphic inverse of $f$ is
\[ (s-x)^{-\prodL} = (x^2 - 2s_0x + |s|^2)^{-1}(\overline{s}-x) = S_L^{-1}(s,x).\]
Similarly, the right slice hyperholomorphic inverse of $f$ is
\[ (s-x)^{-\prodR} = (\overline{s} - x)(x^2 - 2s_0x + |s|^2)^{-1}=S_R^{-1}(s,x).\]
\end{Ex}

The S-functional calculus was constructed using the notion of slice hyperholomorphic functions given in Definition \ref{defsmon}.
The most important properties of slice hyperholomorphic functions are based on the Representation Formula or Structure Formula which is not easy to deduce from this definition. Such formula is deduced using the hypotheses that the domain on which $f$ is defined
is axially symmetric and it is an s-domain.
From the representation formula it follow that a slice hyperholomorphic function is of the form $f=\alpha+I\beta$. Precisely Corollory 2.2.20 in \cite{css_book} claims:
{\it
Let $U\subseteq \mathbb{R}^{n+1}$  be an axially symmetric s-domain,  and $D\subseteq\rr^2$ be such that $u+Iv\in U$ whenever $(u,v)\in D$ and let $f : U \to \mathbb{R}_{n}$. The function $f$ is a slice monogenic function if and only if
there exist two differentiable functions $\alpha, \beta: D\subseteq\rr^2 \to \mathbb{R}_n$ satisfying $\alpha(u,v)=\alpha(u,-v)$, $\beta(u,v)=-\beta(u,-v)$ and the Cauchy-Riemann system
 \begin{equation}\label{CRSIST}
\left\{
\begin{array}{c}
\pp_u \alpha-\pp_v\beta=0\\
\pp_u \beta+\pp_v\alpha=0\\
\end{array}
\right.
\end{equation}
and such that
\begin{equation}
f(u+Iv)=\alpha(u,v)+I\beta(u,v).
\end{equation}
}
\begin{Rk}\label{RKTWODEF}
{\rm A different approach to define slice hyperholomorphic functions  is to consider functions of the form $f=\alpha+I\beta$ such that
$\alpha$ and $\beta$ satisfying $\alpha(u,v)=\alpha(u,-v)$, $\beta(u,v)=-\beta(u,-v)$ and the Cauchy-Riemann system. In this case the Representation formula, and as a consequence the Cauchy formula can be stated removing the hypothesis that the domain of definition $U$ is an s-domain.
It is enough to require that it is axially symmetric. This is because with this definition the representation formula now is easily obtained
from the definition $f(u+Iv)=\alpha(u,v)+I\beta(u,v)$.

In a sense we can claim that Definition \ref{defsmon} contains more functions then the functions defined as  $f(u+Iv)=\alpha(u,v)+I\beta(u,v)$. To obtain them we need an additional condition on the domain, to be an s-domain, but if we already assume that $f(u+Iv)=\alpha(u,v)+I\beta(u,v)$
plus $\alpha(u,v)=\alpha(u,-v)$, $\beta(u,v)=-\beta(u,-v)$ and the Cauchy-Riemann system then we can remove the condition s-domain.
This approach as been used in \cite{GP} to study slice regularity on real alternative algebras.
}
\end{Rk}
For the purpose of this paper, in order to remove any restriction we need domains that are axially symmetric only
so the slice hyperholomorphicity is intended in this second definition.
In this paper we leave both definitions because the S-functional calculus and its properties are proved with Definition \ref{defsmon} on axially symmetric s-domains, but with the definition in the above remark we can remove the fact that we take s-domains.
For more comments on the various definitions of slice hyperholomorphic functions see  the paper
\cite{Global}.

\section{ The S-functional calculus of non commuting operators}

For the results in this section see the book \cite{css_book} and \cite{acgs}.
In the sequel, we will consider a Banach space $V$ over
$\mathbb{R}$
 with norm $\|\cdot \|$.
It is possible to endow $V$
with an operation of multiplication by elements of $\rr_n$ which gives
a two-sided module over $\rr_n$.
A two-sided module $V$ over $\rr_n$ is called a Banach module over $\rr_n$,
 if there exists a constant $C \geq 1$  such
that $\|va\|\leq C\| v\| |a|$ and $\|av\|\leq C |a|\| v\|$ for all
$v\in V$ and $a\in\rr_n$.
 By $V_n$ we denote $V\otimes \rr_n$ over $\rr_n$; $V_n$ turns out to be a
 two-sided Banach module.\\
 An element in $V_n$ is of the type $\sum_A v_A\otimes e_A$ (where
 $A=i_1\ldots i_r$, $i_\ell\in \{1,2,\ldots, n\}$, $i_1<\ldots <i_r$ is a multi-index).
The multiplications of an element $v\in V_n$ with a scalar
$a\in \rr_n$ are defined by $va=\sum_A v_A \otimes (e_A a)$ and $av=\sum_A v_A \otimes (ae_A )$.
For simplicity, we will write
$\sum_A v_A e_A$ instead of $\sum_A v_A \otimes e_A$. Finally, we define $\| v\|^2_{V_n}=
\sum_A\| v_A\|^2_V$.

We denote by
$\mathcal{B}(V)$  the space
of bounded $\mathbb{R}$-homomorphisms of the Banach space $V$ to itself
 endowed with the natural norm denoted by $\|\cdot\|_{\mathcal{B}(V)}$.
Given $T_A\in \mathcal{B}(V)$, we can introduce the operator $T=\sum_A T_Ae_A$ and
its action on $v=\sum v_Be_B\in V_n$ as $T(v)=\sum_{A,B}
T_A(v_B)e_Ae_B$. The operator $T$ is a right-module homomorphism which is a bounded linear
map on $V_n$.
\\
In the sequel, we will consider operators of the form
$T=T_0+\sum_{j=1}^ne_jT_j$ where ${T_j\in\mathcal{B}(V)}$ for $j=0,1,\ldots ,n$.
The subset of such operators in ${\mathcal{B}(V_n)}$ will be denoted by $\mathcal{B}^{\small 0,1}(V_n)$.
We define $\|T\|_{\mathcal{B}^{\small 0,1}(V_n)}=\sum_j \|T_j\|_{\mathcal{B}(V)}$.
Note that, in the sequel, we will omit the subscript $\mathcal{B}^{\small 0,1}(V_n)$ in the norm of an operator. Note also that  $\|TS\|\leq \|T\| \|S\|$.

\begin{Tm}\label{Ssinistro}
Let $T\in\mathcal{B}^{\small 0,1}(V_n)$ and let $s \in \rr^{n+1}$.
Then,  for $\|T\|< |s|$, we have
\begin{equation}\label{SERLEFTgf}
\sum_{n=0}^\infty T^n s^{-1-n}=-(T^2-2\Re(s) T+|s|^2\mathcal{I})^{-1}(T-\overline{s}\mathcal{I}),
\end{equation}
\begin{equation}\label{SERRIGHTgf}
\sum_{n=0}^\infty s^{-1-n} T^n =-(T-\overline{s}\mathcal{I})(T^2-2\Re(s) T+|s|^2\mathcal{I})^{-1}.
\end{equation}
\end{Tm}
We observe that the fact that these formulas hold is independent of whether the components of the paravector operator $T$ commute or not. Moreover the operators on the right hand sides of
(\ref{SERLEFTgf}) and (\ref{SERRIGHTgf}) are defined on a subset of $\mathbb{R}^{n+1}$
that is larger than $\{ s\in \mathbb{R}^{n+1} \ :\ \|T\|< |s|\}.$ This fact suggests the definition of the $S$-spectrum, of the $S$-resolvent set and of the $S$-resolvent operators.

\begin{Dn}[The $S$-spectrum and the $S$-resolvent set]
Let $T\in\mathcal{B}^{\small 0,1}(V_n)$.
We define the $S$-spectrum $\sigma_S(T)$ of $T$  as
$$
\sigma_S(T):=\{ s\in \mathbb{R}^{n+1}\ \ :\ \ T^2-2 \Re(s)T+|s|^2\mathcal{I}\ \ \
{\it is\ not\  invertible}\}.
$$
The $S$-resolvent set $\rho_S(T)$ is defined as
$$
\rho_S(T):=\mathbb{R}^{n+1}\setminus\sigma_S(T).
$$
\end{Dn}

\begin{Dn}[The $S$-resolvent operators]
For  $T\in \mathcal{B}^{\small 0,1}(V_n)$ and $s\in \rho_S(T)$, we define the left $S$-resolvent operator as
\begin{equation}\label{quatSresolrddlft}
S_L^{-1}(s,T):=-(T^2-2\Re(s) T+|s|^2\mathcal{I})^{-1}(T-\overline{s}\mathcal{I}),
\end{equation}
and the right $S$-resolvent operator as
\begin{equation}\label{quatSresorig}
S_R^{-1}(s,T):=-(T-\overline{s}\mathcal{I})(T^2-2\Re(s) T+|s|^2\mathcal{I})^{-1}.
\end{equation}
\end{Dn}
The operators $S_L^{-1}(s,T)$ and $S_R^{-1}(s,T)$ satisfy the equations below, see \cite{acgs}.
\begin{Tm}\label{RLRESOLVEQ}
Let $T\in\mathcal{B}^{\small 0,1}(V_n)$ and let $s,p\in \rho_S(T)$. Then we have
\begin{equation}\label{RLresolv}\small
S_R^{-1}(s,T)S_L^{-1}(p,T)=\big[(S_R^{-1}(s,T)-S_L^{-1}(p,T))p
-\overline{s}(S_R^{-1}(s,T)-S_L^{-1}(p,T))\big](p^2-2s_0p+|s|^2)^{-1}.
\end{equation}
Moreover, the resolvent equation can also be written as
\begin{equation}\label{RLresolvII}\small
S_R^{-1}(s,T)S_L^{-1}(p,T)=(s^2-2p_0s+|p|^2)^{-1}\big[s(S_R^{-1}(s,T)-S_L^{-1}(p,T))
-(S_R^{-1}(s,T)-S_L^{-1}(p,T))\overline{p} \big].
\end{equation}
\end{Tm}

Here we recall the formulations of the S-functional calculus.
We first recall two important properties of
the S-spectrum.
\begin{Tm}[Structure of the $S$-spectrum]\label{strutturaS}
\par\noindent
Let $T\in\mathcal{B}^{\small 0,1}(V_n)$
and suppose that $p=p_0+\underline{p}$ belongs  $\sigma_S(T)$
 with $ \underline{p}\neq 0$.
Then all the elements of the $(n-1)$-sphere $[p]$ belong to $\sigma_S(T)$.
\end{Tm}
This result implies that if $p\in\sigma_S(T)$ then either $p$ is a
real point or the whole $(n-1)$-sphere $[p]$ belongs to
$\sigma_S(T)$.

\par\noindent
\begin{Tm}[Compactness of the $S$-spectrum]\label{compattezaS}
 Let $T\in\mathcal{B}^{\small 0,1}(V_n)$. Then
 $\sigma_S (T)$  is a compact nonempty set that is
contained in the closed ball $\{{s\in\rr^{n+1} :\,}  {|s|\leq \|T\|} \}$.
\end{Tm}

\begin{Dn}
  Let $V_n$ be a two-sided Banach module, let  $T\in\mathcal{B}^{0,1}(V_n)$and let $U \subset \mathbb{R}^{n+1}$ be an axially symmetric s-domain
that contains  the $S$-spectrum $\sigma_S(T)$  such that
$\partial (U\cap \mathbb{C}_I)$ is the union of a finite number of
continuously differentiable Jordan curves  for every $I\in\mathbb{S}$.
In this case we say that $U$ is a $T$-admissible domain.
\end{Dn}
We can now introduce the class of functions for which we can define the two versions of the S-functional calculus.
 \begin{Dn}\label{quatdef3.9}
Let $V_n$ be a two-sided Banach module and let  $T\in\mathcal{B}^{0,1}(V_n)$.
\begin{enumerate}[(i)]
\item
A function  $f$  is said to be locally left slice hyperholomorphic  on $\sigma_S(T)$
if there exists a $T$-admissible domain $U\subset \mathbb{R}^{n+1}$ such that $f$ is left slice hyperholomorphic on some open set $W$ with $\overline{U}\subset W$.
We denote by $\mathcal{SM}^L_{\sigma_S(T)}$ the set of locally
left slice hyperholomorphic functions on $\sigma_S (T)$.
\item
A function  $f$  is said to be locally right slice hyperholomorphic  on $\sigma_S(T)$
if there exists a $T$-admissible domain $U\subset \mathbb{R}^{n+1}$ such that $f$ is right slice hyperholomorphic on some open set $W$ with $\overline{U}\subset W$.
We denote by $\mathcal{SM}^R_{\sigma_S(T)}$ the set of locally
right  slice hyperholomorphic functions on $\sigma_S (T)$.
\item We define the class $\bimon_{\sigma_S(T)}$ as the set of all locally left slice hyperholomorphic functions on $\sigma_S(T)$ that satisfy $f(U\cap\C_I)\subset\C_I, I\in\S$ for some axially symmetric slice domain $U$ with $\sigma_S(T)\subset U$.
\end{enumerate}
\end{Dn}

\begin{Dn}[The  $S$-functional calculus]
Let $V_n$ be a two-sided  Banach module and  ${T\in\mathcal{B}^{0,1}(V_n)}$.
  Let $U\subset  \mathbb{R}^{n+1}$ be a $T$-admissible domain and set $ds_I=- ds I$. We define
\begin{equation}\label{Scalleft}
f(T)={{1}\over{2\pi }} \int_{\partial (U\cap \mathbb{C}_I)} S_L^{-1} (s,T)\  ds_I \ f(s), \ \ {\it for}\ \ f\in \mathcal{SM}^L_{\sigma_S(T)},
\end{equation}
and
\begin{equation}\label{Scalright}
f(T)={{1}\over{2\pi }} \int_{\partial (U\cap \mathbb{C}_I)} \  f(s)\ ds_I \ S_R^{-1} (s,T),\ \  {\it for}\ \ f\in \mathcal{SM}^R_{\sigma_S(T)}.
\end{equation}
\end{Dn}

\begin{La}\label{PropSCalc}
Let $V_n$ be a two-sided  Banach module and  $T\in\mathcal{B}^{0,1}(V_n)$.
\begin{enumerate}[(i)]
\item Let $f,g\in\SM_{\sigma_S(T)}^L$ and let $a\in\H$. Then $(fa+g)(T) = f(T)a + g(T)$. Moreover, if $f\in\bimon_{\sigma_S(T)}$, then $(fg)(T) = f(T)g(T)$.
\item Let $f,g\in\SM_{\sigma_S(T)}^R$ and let $a\in\H$. Then $(af+g)(T) = af(T) + g(T)$. Moreover, if $g\in\bimon_{\sigma_S(T)}$, then $(fg)(T) = f(T)g(T)$.
\item Let $f_n, f \in\SM^L_{\sigma_S(T)}$ or $f_n, f \in\SM^R_{\sigma_S(T)}$. If there exists a $T$-admissible  domain $U$ such that $f_n\to f$ uniformly on $\overline{U}$, then $f_n(T) \to f(T)$ in the operator norm.
\end{enumerate}
\end{La}

\section{Another proof of the invertibility of $T^2-\Re(s)T + |s|^2\mathcal I$ for $\|T\|<|s|$}
To determine the closed form of the noncommutative Cauchy kernel series as in Theorem~3.1.3 of \cite{css_book}, traditionally, one shows the identity
\[ (T^2 - 2 \Re(s) T + |s|^2\id)\sum_{n\geq 0} T^ns^{-1-n} = \overline{s} \id - T.\]
For $\|T\|<|s|$, the invertability of $\overline{s}\:\!\id-T$
implies,
\[
-(T-\overline{s}\:\!\id)^{-1}(T^2 - 2 \Re(s) T + |s|^2\id)\sum_{n\geq 0} T^ns^{-1-n} = \id.\]
The following computations show explicitly that the identity
\begin{equation}\label{MainSeries}
\left(\sum_{n=0}^{\infty}T^ns^{-n-1} \right)(\overline{s}\:\!\id-T)^{-1}(T^2 - 2\Re(s)T + |s|^2\id) =\id\end{equation}
also holds true for $\|T\| < |s|$.

Recall that  by Theorem~4.14.6 in \cite{css_book}, we have:
$$
(\overline{s}\:\!\id-T)^{-1} = \sum_{n=0}^{\infty} (\overline{s}^{-1}T)^n\overline{s}^{-1}\ \  {\rm for}\ \ \ \|T\|<|s|.
$$
Therefore
\begin{align*}\left(\sum_{n=0}^{\infty}T^ns^{-n-1}  \right)&(\overline{s}\:\!\id-T)^{-1}(T^2 - 2\Re(s)T + |s|^2\id) \\
&= \left(\sum_{n=0}^{\infty}T^ns^{-n-1} \right)\left(\sum_{n=0}^{\infty} \bar s^{-1} (T\bar s^{-1})^n\right)(T^2 - 2\Re(s)T + |s|^2\id) \\
&=\sum_{n=0}^{\infty}\sum_{k=0}^{n} T^ks^{-k-1}\bar s^{-1} (T\bar s^{-1})^{n-k} (T^2 - 2\Re(s)T + |s|^2\id) \\
&=\sum_{n=0}^{\infty}\sum_{k=0}^{n} T^ks^{-k} (T\bar s^{-1})^{n-k} \left(\frac{1}{|s|^2}T^2 - (s^{-1}+\bar{s}^{-1})T + \id\right),
\end{align*}
where the last equality follows from $s^{-1}\overline{s}^{-1} = 1/|s|^2$.
To simplify this sum, we need two auxiliary lemmas.

\begin{La} \label{LemmaRPFL} Let $T\in\mathcal{B}^{0,1}(V_n)$. For any $s\in\mathbb{H}$, we have
\begin{equation}\label{RPFL}\sum_{k=0}^{n}T^ks^{-k}(T\overline{s}^{-1})^{n-k} =T^{n}\sum_{k=0}^{n}s^{-k}\overline{s}^{-n+k}.\end{equation}
\end{La}
\begin{proof}
For $n=0$ the identity \eqref{RPFL} is trivial. If we assume that it holds true for $n-1$, then
\[
\begin{split}
\sum_{k=0}^{n}T^ks^{-k}(T\overline{s}^{-1})^{n-k} &= \left(\sum_{k=0}^{n-1}T^ks^{-k}(T\overline{s}^{-1})^{n-1-k}\right) T\overline{s}^{-1} + T^ns^{-n}
\\
&
= \left(T^{n-1}\sum_{k=0}^{n-1}s^{-k}\overline{s}^{-(n-1)+k}\right)T\overline{s}^{-1} + T^n s^{-n}
\\
&
= T^{n}\sum_{k=0}^{n-1}s^{-k}\overline{s}^{-n+k} + T^n s^{-n}
\\
&
 = T^{n}\sum_{k=0}^{n}s^{-k}\overline{s}^{-n+k},
\end{split}
\]
where we used that $\sum_{k=0}^{n-1}s^{-k}\overline{s}^{-(n-1)+k}$ is real and therefore commutes with $T$ because
\[\overline{\sum_{k=0}^{n-1}s^{-k}\overline{s}^{-(n-1)+k}} = \sum_{k=0}^{n-1}s^{-(n-1)+k}\overline{s}^{k} = \sum_{k=0}^{n-1}s^{-k}\overline{s}^{-(n-1)+k}.\]
\end{proof}

\begin{La}\label{RPFL2}
Let $T\in\mathcal{B}^{0,1}(V_n)$ and $s\in\mathbb{H}$. If we define
\begin{equation}\label{S(N)}S(m) = \sum_{n=0}^{m}\sum_{k=0}^{n} T^ks^{-k} (T\bar s^{-1})^{n-k} \left(\frac{1}{|s|^2}T^2 - (s^{-1}+\bar{s}^{-1})T + \id\right)
\end{equation}
for $m\in\N_0$, then
\begin{equation}\label{S(N)id}S(m) = \id - T^{m+1}\sum_{k=0}^{m+1}s^{-k}\overline{s}^{-(m+1)+k} + T^{m+2}\frac{1}{|s|^2}\sum_{k=0}^{m}s^{-k}\overline{s}^{-m+k}.\end{equation}
\end{La}
\begin{proof}
From Lemma \ref{LemmaRPFL}, we deduce

\begin{align*}S(m) =& \sum_{n=0}^{m}\sum_{k=0}^{n} T^ks^{-k} (T\bar s^{-1})^{n-k} \left(\frac{1}{|s|^2}T^2 - (s^{-1}+\bar{s}^{-1})T + \id\right) \\
=& \sum_{n=0}^{m}T^n\sum_{k=0}^{n}s^{-k} \overline{s}^{-n+k} \left(\frac{1}{|s|^2}T^2 - (s^{-1}+\bar{s}^{-1})T + \id\right).
\end{align*}
Since $\sum_{k=0}^{n}s^{-k} \overline{s}^{-n+k}$ is real, it commutes with $T$. Together with the identity ${1/|s|^2} = s^{-1}\overline{s}^{-1}$, this implies
\begin{align*}
S(m)=& \sum_{n=0}^m T^{n+2}  \frac{1}{|s|^2} \sum_{k=0}^{n}s^{-k}\overline{s}^{-n+k} - \sum_{n=0}^m T^{n+1}(s^{-1}+\overline{s}^{-1}) \sum_{k=0}^{n}s^{-k} \overline{s}^{-n+k} + \sum_{n=0}^mT^n \sum_{k=0}^{n}s^{-k} \overline{s}^{-n+k} \\
=& \sum_{n=0}^m T^{n+2} \sum_{k=0}^{n}s^{-(k+1)}\overline{s}^{-(n+1) +k} - \sum_{n=0}^m T^{n+1}\sum_{k=0}^{n}s^{-(k+1)} \overline{s}^{-n+k}  \\
&\phantom{\sum_{n=0}^m T^{n+2} \sum_{k=0}^{n}s^{-(k+1)}\overline{s}^{-(n+1) +k}}-  \sum_{n=0}^m T^{n+1} \sum_{k=0}^{n}s^{-k} \overline{s}^{-(n+1)+k} + \sum_{n=0}^mT^n \sum_{k=0}^{n}s^{-k} \overline{s}^{-n+k}.
\end{align*}
By shifting the index $n$ and then the index $k$, we obtain
\begin{align*}
S(m)=& \sum_{n=2}^{m+2} T^{n} \sum_{k=0}^{n-2}s^{-(k+1)}\overline{s}^{-(n-1)+k} - \sum_{n=1}^{m+1} T^{n}\sum_{k=0}^{n-1}s^{-(k+1)} \overline{s}^{-(n-1)+k}  \\
&-  \sum_{n=1}^{m+1} T^{n} \sum_{k=0}^{n-1}s^{-k} \overline{s}^{-n+k} + \sum_{n=0}^mT^n \sum_{k=0}^{n}s^{-k} \overline{s}^{-n+k}\\
=& \sum_{n=2}^{m+2} T^{n} \sum_{k=1}^{n-1}s^{-k}\overline{s}^{-n+k} - \sum_{n=1}^{m+1} T^{n}\sum_{k=1}^{n}s^{-k} \overline{s}^{-n+k}  \\
&-  \sum_{n=1}^{m+1} T^{n} \sum_{k=0}^{n-1}s^{-k} \overline{s}^{-n+k} + \sum_{n=0}^mT^n \sum_{k=0}^{n}s^{-k} \overline{s}^{-n+k}.
\end{align*}
Adding the terms with index $k=n$ of the second sum to the third sum yields
\begin{align*}
S(m)=& \sum_{n=2}^{m+2} T^{n} \sum_{k=1}^{n-1}s^{-k}\overline{s}^{-n+k} - \sum_{n=2}^{m+1} T^{n}\sum_{k=1}^{n-1}s^{-k} \overline{s}^{-n+k}  \\
&-  \sum_{n=1}^{m+1} T^{n} \sum_{k=0}^{n}s^{-k} \overline{s}^{-n+k} + \sum_{n=0}^mT^n \sum_{k=0}^{n}s^{-k} \overline{s}^{-n+k}\\
=& T^{m+2}\sum_{k=1}^{m+1}s^{-k}\overline{s}^{-(m+2)+k} - T^{m+1}\sum_{k=0}^{m+1}s^{-k}\overline{s}^{-(m+1)+k} + \id \\
=&T^{m+2}\frac{1}{|s|^2}\sum_{k=0}^{m}s^{-k}\overline{s}^{-m+k} - T^{m+1}\sum_{k=0}^{m+1}s^{-k}\overline{s}^{-(m+1)+k} + \id .
\end{align*}
\end{proof}

\begin{La}
Let $T\in\mathcal{B}^{0,1}(V_n)$ and $s\in\mathbb{H}$ with $\|T\|  < |s|$. Then the identity \eqref{MainSeries} holds true, that is,
\[\left(\sum_{n=0}^{\infty}T^ns^{-n-1} \right)(\overline{s}\:\!\id-T)^{-1}(T^2 - 2\Re(s)T + |s|^2\id) =\id.\]
\end{La}
\begin{proof}
The considerations before Lemma~\ref{LemmaRPFL} showed
\[
\begin{split}
&
\left(\sum_{n=0}^{\infty}T^ns^{-n-1}  \right)(\overline{s}\:\!\id-T)^{-1}(T^2 - 2\Re(s)T + |s|^2\id)
\\
&
 =\sum_{n=0}^{\infty}\sum_{k=0}^{n} T^ks^{-k} (T\bar s^{-1})^{n-k} \left(\frac{1}{|s|^2}T^2 - (s^{-1}+\bar{s}^{-1})T + \id\right).
 \end{split}
 \]
Lemma~\ref{RPFL2} implies
\[
\begin{split}
&\left\|\left(\sum_{n=0}^{\infty}T^ns^{-n-1} \right)(\overline{s}\:\!\id-T)^{-1}(T^2 - 2\Re(s)T + |s|^2\id) - \id\right\|
 \\
 &= \lim_{m\to\infty} \|S(m) - \id\|
 \\
&
= \lim_{m\to\infty}\left\| - T^{m+1}\sum_{k=0}^{m+1}s^{-k}\overline{s}^{-(m+1)+k} + T^{m+2}\frac{1}{|s|^2}\sum_{k=0}^{m}s^{-k}\overline{s}^{-m+k}\right\|
\\
&
\leq \lim_{m\to\infty} \|T\|^{m+1}\sum_{k=0}^{m+1}|s|^{-k}|\overline{s}|^{-(m+1)+k} + \|T\|^{m+2}\frac{1}{|s|^2}\sum_{k=0}^{m}|s|^{-k}|\overline{s}|^{-m+k}
\\
&
= \lim_{m\to\infty} (m+1)\frac{\|T\|^{m+1}}{|s|^{m+1}} + (m+2)\frac{\|T\|^{m+2}}{|s|^{m+2}} = 0
 \end{split}
 \]
because $\|T\| < |s|$. Hence, \eqref{MainSeries} holds true.

\end{proof}

\section{Slice hyperholomorphic powers of the resolvent operator}
Recall that the product rule $(fg)(T) = f(T)g(T)$ for the $S$-functional calculus does only hold true if $f\in\bimon_{\sigma_S(T)}$ and $g\in\SM_{\sigma_S(T)}^L$ or if $f\in\SM^R_{\sigma_S(T)}$ and $g\in\bimon_{\sigma_S(T)}$. This is due to the fact that, for $f,g\in\SM_{\sigma_S(T)}^L$ or for $f,g\in\SM_{\sigma_S(T)}^R$, the product $fg$ does in general not belong to $\SM_{\sigma_S(T)}^L$ resp. $\SM_{\sigma_S(T)}^R$. If on the other hand one considers the left slice hyperholomorphic product $f\prodL g$ of two left slice hyperholomorphic functions, then it is not clear to which operation between operators it corresponds. Similarly, in the right slice hyperholomorphic case, it is not clear how to define the $\prodR$-product between two operators in general. However, at least for power series of an operator variable, we can use the formulas \eqref{prodLSeries} and \eqref{prodRSeries} to define a their $\prodL$- resp. $\prodR$-product.
\begin{Dn}\label{OpSliceProd}
Let $T\in\boundOP^{0,1}(V_n)$.
For $F = \sum_{n=0}^\infty T^na_n$
and
$G = \sum_{n=0}^\infty T^nb_n$
with ${a_\ell,b_\ell\in\R_{n}}$ for $\ell\in \mathbb{N}$,
we define
\[
F\prodL G := \sum_{n=0}^\infty  T^n\left(\sum_{k=0}^na_{k}b_{n-k}\right).
 \]
For $\widetilde{F} = \sum_{n=0}^\infty a_nT^n$ and $\widetilde{G} = \sum_{n=0}^\infty b_nT^n$, we define
\[ \widetilde{F}\prodR \widetilde{G} := \sum_{n=0}^\infty  \left(\sum_{k=0}^na_{k}b_{n-k}\right)T^n. \]
\end{Dn}
\begin{Rk}\label{SeriesProdSimple}
Note that $F\prodL G = FG$ if $a_n\in\R$ for any $n\in\N$. In this case, the coefficients $a_n$ commute with the operator $T$, and hence,
\[ F\prodL G = \sum_{n=0}^\infty  T^n\left(\sum_{k=0}^na_{k}b_{n-k}\right) = \sum_{n=0}^{\infty}\sum_{k=0}^nT^ka_k T^{n-k}b_{n-k} = F G.\]
Similarly, $\widetilde{F}\prodR\widetilde{G} = \widetilde{F}\widetilde{G}$ if $b_n\in\R$ for any $n\in\N$.
\end{Rk}

\begin{Cy}\label{SliceProdRule}
Let $T\in\boundOP^{0,1}(V_n)$ and let $f(x) = \sum_{n=0}^\infty x^na_n$ and $g(x) = \sum_{n=0}^\infty x^nb_n$ be two left slice hyperholomorphic power series that converge on a ball $B_r(0)$, of radius $r>0$ and centered at $0$, with $\sigma_S(T)\subset B_r(0)$. Then
\[ f(T)\prodL g(T)  = (f\prodL g)(T).\]
Similarly, for two right slice hyperholomorphic power series $\tilde{f}(x) = \sum_{n=0}^\infty a_nx^n$ and $\tilde{g}(x) = \sum_{n=0}^\infty b_nx^n$  that converge on a ball $B_r(0)$ with $\sigma_S(T)\subset B_r(0)$, we have
\[ \tilde{f}(T)\prodR \tilde{g}(T) = (\tilde{f}\prodR \tilde{g})(T).\]
\end{Cy}
\begin{proof}
By the properties of the $S$-functional calculus, we have $f(T) = \sum_{n=0}^\infty T^na_n$ and $g(T) = \sum_{n=0}^\infty T^nb_n$. Hence,
\begin{align*}
f(T)\prodL g(T) &= \sum_{n=0}^\infty T^n\left(\sum_{k=0}^n a_k b_{n-k}\right) \\
&= \frac{1}{2\pi}\int_{\partial(U\cap\C_I)} S_L^{-1}(s,T)\,ds_I\, \sum_{n=0}^\infty s^n\left(\sum_{k=0}^n a_k b_{n-k}\right)\\
&= \frac{1}{2\pi}\int_{\partial(U\cap\C_I)}S_L^{-1}(s,T)\,ds_I\, f\prodL g(s) = (f\prodL g)(T).
\end{align*}
An analogous computation shows the right slice hyperholomorphic case.

\end{proof}

Recall from Example~\ref{ExSliceInv} that $S_L^{-1}(s,x) = (s-x)^{-\prodL}$ and that $S_R^{-1}(s,x) = (s-x)^{-\prodR}$. This motivates the following definition.
\begin{Dn}
Let $s\in\H$. For $n\in\Z$, we define
\[ S_L^n(s,x) := (s-x)^{\prodL n}\qquad and \qquad S_R^n(s,x) :=(s-x)^{\prodR n}.\]

\end{Dn}

\begin{Cy}\label{ScalarPowers}
Let $s\in\H$. For $n\geq 0$, we have
\begin{equation}\label{PosPow} S_L^n (s,x) = \sum_{k=0}^n\binom{n}{k} (-x)^k s^{n-k}\qquad\text{and}\qquad S_R^n(s,x) = \sum_{k=0}^n\binom{n}{k} s^{n-k}(-x)^k\end{equation}
and
\begin{equation}\label{NegPow}S_L^{-n}(s,x) = (x^2 - 2s_0 x + |s|^2)^{-n} (\overline{s}-x)^{\prodL n}
 \quad\text{and}\quad S_R^{-n}(s,x) = (\overline{s}-x)^{\prodR n}(x^2 - 2s_0 x + |s|^2)^{-n}.\end{equation}
Moreover, for $m,n\geq 0$, we have
\[ S_L^{-n}(s,x)\prodL S_L^{-m}(\overline{s},x) = (x^2 - 2s_0 x + |s|^2)^{-(n+m)}\left[(\overline{s}-x)^{\prodL n}\prodL({s}-x)^{\prodL m}\right]\]
and
\[ S_R^{-n}(s,x)\prodR S_R^{-m}(\overline{s},x) = \left[(\overline{s}-x)^{\prodR n}\prodR(s-x)^{\prodR m}\right](x^2 - 2s_0 x + |s|^2)^{-(n+m)}.\]
\end{Cy}
\begin{proof}
For $n=0$, we have $(s-x)^{\prodL 0} = 1$, and hence \eqref{PosPow} is obviously true. Assume that it holds true for $n-1$. Then \eqref{prodLSeries} implies
\begin{align*} S_L^n(s,x) &= (s-x)^{\prodL n} = (s-x)^{\prodL (n-1)}\prodL (s-x)  =  (s-x)^{\prodL (n-1)}\prodL s+  (s-x)^{\prodL (n-1)}\prodL (-x)  \\
&= \sum_{k=0}^{n-1}\binom{n-1}{k} (-x)^ks^{n-k} + \sum_{k=0}^{n-1}\binom{n-1}{k}(-x)^{k+1}s^{n-1-k} \\
&= \sum_{k=0}^{n-1}\binom{n-1}{k} (-x)^ks^{n-k} + \sum_{k=1}^{n}\binom{n-1}{k-1}(-x)^{k}s^{n-k} = \sum_{k=0}^n\binom{n}{k}(-x)^ks^{n-k}
\end{align*}
and \eqref{PosPow} follows by induction.

We also show the identity \eqref{NegPow} by induction. It is obviously true for $n=0$. Assume that it holds true for $n-1$ and observe that $(x^2 - 2s_0x + |s|^2)^{-1}\in\bimon(\H\setminus{[s]})$.  Corollary~\ref{BimonProp}, (\ref{BimonProd}) yields
\begin{align*} S_L^{-n}(s,x)& = (s-x)^{-\prodL (n-1)} \prodL (s-x)^{-\prodL}\\
& = \left[(x^2 - 2s_0 x + |s|^2)^{-(n-1)}(\overline{s}-x)^{\prodL (n-1)}\right]\prodL\left[ (x^2 - 2s_0x + |s|^2)^{-1}(\overline{s}-x)\right] \\
& = (x^2 - 2s_0x + |s|^2)^{-(n-1)}\prodL(\overline{s}-x)^{\prodL (n-1)}\prodL(x^2-2s_0x + |s|^2)^{-1}\prodL (\overline{s}-x) \\
&= (x^2 - 2s_0x + |s|^2)^{-(n-1)}\prodL(x^2-2s_0x + |s|^2)^{-1}\prodL(\overline{s}-x)^{\prodL (n-1)}\prodL (\overline{s}-x)  \\
&= (x^2 - 2s_0x + |s|^2)^{-n}(\overline{s}-x)^{\prodL n}.
\end{align*}
Finally, (\ref{BimonProd}) in  Corollary~\ref{BimonProp} also implies  for $m,n\geq 0 $
\begin{align*} S_L^{-n}(s,x) \prodL S_L^{-m}(\overline{s},x) &= \left[(x^2 - 2s_0x + |s|^2)^{-n}(\overline{s}-x)^{\prodL n}\right]\prodL \left[(x^2 - 2s_0x + |s|^2)^{-m}({s}-x)^{\prodL m}\right]\\
&= (x^2 - 2s_0x + |s|^2)^{-n}\prodL(\overline{s}-x)^{\prodL n}\prodL (x^2 - 2s_0x + |s|^2)^{-m}\prodL({s}-x)^{\prodL m}\\
&= (x^2 - 2s_0x + |s|^2)^{-n}\prodL (x^2 - 2s_0x + |s|^2)^{-m}\prodL (\overline{s}-x)^{\prodL n}\prodL ({s}-x)^{\prodL m}\\
&= (x^2 - 2s_0x + |s|^2)^{-(n+m)}\left[(\overline{s}-x)^{\prodL( n+m)}\prodL(s-x)^{\prodL m}\right].
\end{align*}
The right slice hyperholomorphic case can be shown by similar computations.

\end{proof}

\begin{Dn}\label{ResSlicePowers}
Let $T\in\boundOP^{0,1}(V_n)$ and let $s\in\rho_S(T)$. For $n,m\geq 0$, we define
\[S_L^{-n}(s,T):= (T^2 - 2s_0T + |s|^2\id)^{-n}(\overline{s}\:\!\id-T)^{\prodL n}\]
and
\[S_L^{-n}(s,T)\prodL S_L^{-m}(\overline{s},T):= (T^2 - 2s_0T + |s|^2\id)^{-(n+m)}\left[(\overline{s}\:\!\id - T)^{\prodL n}\prodL(s\id-T)^{\prodL m}\right].\]
Similarly, we define
\[S_R^{-n}(s,T):= (\overline{s}\id-T)^{\prodR n} (T^2 - 2s_0T + |s|^2\id)^{-n}\]
and
\[S_R^{-n}(s,T)\prodR S_R^{-m}(\overline{s},T):= \left[(\overline{s}\:\!\id - T)^{\prodR n}\prodR(s\id-T)^{\prodR m}\right](T^2 - 2s_0T + |s|^2\id)^{-(n+m)}.\]
\end{Dn}
\begin{Rk}
Because of Lemma~\ref{ScalarPowers} and Corollary~\ref{PropSCalc}, this definition is consistent with the $S$-functional calculus in the sense that
$$
S_L^{-n}(s,T) = \left[S_L^{-n}(s,\cdot)\right](T)
$$
and
$$
S_L^{-n}(s,T)\prodL S_L^{-m}(\overline{s},T) = \left[ S_L^{-n}(s,\cdot)\prodL S_L^{-m}(\overline{s},\cdot)\right](T).
 $$
 However, note that the notations $S_L^{-n}(s,T)$ and $S_L^{-n}(s,T)\prodL S_L^{-m}(\overline{s},T)$ are purely formal. They can only be interpreted as the $\prodL$-products of the operators $S_L^{-n}(s,T)$ and $S_L^{-m}(\overline{s},T)$, if $|s|$ is greater than the spectral radius $r(T)=\sup\{|s|:s\in\sigma_S(T)\}$ of $T$. Otherwise, the operators $S_L^{-1}(s,T)$ and $S_L^{-1}(\overline{s},T)$ do not admit a power series representation and the respective $\prodL$-products  are  not defined.

\end{Rk}

We conclude this section with a corollary that will be useful in several proofs of the following section.
\begin{Cy}\label{PowersSimple}
Let $s=s_0 + I s_1\in\H$ and $n,m\in\N_0$. If $x\in\C_I$, then
\begin{equation}\label{PosPowSimple}(s-x)^{\prodL m}\prodL (\overline{s} - x)^{\prodL n} = (s-x)^m(\overline{s}-x)^n\end{equation}
and
\[S_L^{-m}(s,x)\prodL S_L^{-n}(\overline{s},x) =  (s-x)^{-m}(\overline{s}-x)^{-n}.\]
Moreover, for any $n\in\N_0$, the function
\begin{equation}\label{PowersSumScalar}\sum_{k=0}^n (\overline{s}-x)^{\prodL(k+1)}\prodL(s-x)^{\prodL(n-k+1)}\end{equation}
is a polynomial with real coefficients.
\end{Cy}
\begin{proof}
If $x\in\C_I$, then $s$, $\overline{s}$ and $x$ commute. Hence, it follows from \eqref{PosPow} and the binomial theorem that $(s-x)^{\prodL m} = (s-x)^m$ and $(\overline{s}-x)^{\prodL n} = (\overline{s}-x)^n$. From Remark~\ref{RemarkSimpleProduct}, we deduce that $\eqref{PosPowSimple}$ holds true.
Since $x$ and $s$ commute, we also have
\[\begin{split}S_L^{-m}(s,x) &= (x^2 - 2s_0x + |s|^2)^{-m}(\overline{s}-x)^{\prodL m} \\
&= \left((s-x)(\overline{s}-x)\right)^{-m}\sum_{k=0}^m\binom{m}{k}(-x)^{k}\overline{s}^{m-k} \\
&= (s-x)^{-m}(\overline{s}-x)^{-m} (\overline{s}-x)^m = (s-x)^{-m}.\end{split}\]
An analogous computation shows $S_L^{-n}(\overline{s},x) = (\overline{s}-x)^{-n}$.  Remark~\ref{RemarkSimpleProduct} implies
\[S_L^{-m}(s,x)\prodL S_L^{-n}(\overline{s},x) =  (s-x)^{-m}\prodL(\overline{s}-x)^{-n} =  (s-x)^{-m}(\overline{s}-x)^{-n}.\]
Finally, we consider the sum \eqref{PowersSumScalar}. We denote
\[P(x) := \sum_{k=0}^n (\overline{s}-x)^{\prodL(k+1)}\prodL(s-x)^{\prodL(n-k+1)}.\] The restriction $P_I$ of this function to the plane $\C_I$ is the complex polynomial $P_I(x) = \sum_{k=0}^n(\overline{s}-x)^{k+1}(s-x)^{n-k+1}$. From  the relation
\[ P_I(\overline{x}) = \sum_{k=0}^n(\overline{s}-\overline{x})^{k+1}(s-\overline{x})^{n-k+1} = \overline{\sum_{k=0}^n(s-{x})^{k+1}(\overline{s}-{x})^{n-k+1}}=\overline{P_I(x)},\]
we deduce that its coefficients are real. Consequently, $P = \ext_L(P_I)$ is a polynomial with real coefficients on $\R^{n+1}$.

\end{proof}

\section{Power series expansions in the operator variable}
The following result clarifies how one has to measure the distance between a point $s\in\rho_S(T)$ and the $S$-spectrum of $T$.
\begin{La} \label{DistLem}Let $A\subset\rr^{n+1}$ be axially symmetric and  let $s = s_0 + Is_1\in\H$. Then
\[\dist(s,A) =  \dist(s,A\cap\C_I),\]
where $\dist(s,A) := \inf\{|s-x|:x\in A\}$ denotes the distance between $s$ and the set $A$.
\end{La}
\begin{proof}
 For $x=x_0 + I_x x_1 \in A$ denote $x_I = x_0 + Ix_1$ and chose $J\in\S$ with $I\perp J$ such that $I_x \in \lin\{I, J\}$. Then $x = x_0 + \widetilde{x}_1 I + \widetilde{x}_2 J$ with $\widetilde{x}_1^2 + \widetilde{x}_2^2 = |\underline{x}|^2 = x_1^2$, and in turn
 \begin{align*}
  |s - x_I|^2 &= (s_0-x_0)^2 + (s_1 - x_1)^2 = (s_0-x_0)^2 + s_1^2 - 2s_1x_1 + x_1^2 \\
 &= (s_0-x_0)^2 + s_1^2 - 2 s_1\sqrt{\widetilde{x}_1^2 + \widetilde{x}_2^2} + \widetilde{x}_1^2 + \widetilde{x}_2^2 \\
 &\leq (s_0-x_0)^2 + s_1^2 - 2 s_1\widetilde{x}_1 + \widetilde{x}_1^2 + \widetilde{x}_2^2 = (s_0-x_0)^2 + (s_1-\widetilde{x}_1)^2 + \widetilde{x}_2^2 = |s - x|^2.
 \end{align*}
Since $A$ is axially symmetric, we have $\{x_I: x\in A\}\subset A\cap\C_I$. Consequently,
\[\inf_{x\in A} |s-x| \leq \inf_{x\in A\cap\C_I}|s-x| \leq \inf_{x\in A}|s-x_I| \leq \inf_{x\in A}|s-x|,\]
and in turn,
\[ \dist(s,A) = \inf_{x\in A}|s-x|=  \inf_{x\in A\cap\C_I}|s-x| = \dist(s,A\cap\C_I).\]
\end{proof}

\begin{Pn}\label{KTEstimate}
Let $T\in \mathcal{B}^{0,1}(V_n)$ and let $C\subset\rr^{n+1}$ with $\dist(C,\sigma_S(T))>\varepsilon$ for some $\varepsilon > 0$. Then there exists  a positive constant $K_T$ such that
\begin{equation}\label{sESTSRES}
\left\| S_L^{-m}(s,T)\prodL S_L^{-n}(\overline{s},T)\right\|\leq \frac{K_T}{\varepsilon^{m+n}}
\quad\text{and}\quad
\left\| S_R^{-m}(s,T)\prodL S_R^{-n}(\overline{s},T)\right\|\leq \frac{K_T}{\varepsilon^{m+n}},
\end{equation}
for any $s\in C$ and any $m,n\geq 0$.
\end{Pn}
\begin{proof}
Let $U$ be a $T$-admissible domain, such that $\dist(C, \overline{U}) >\varepsilon$.
We chose $s = {s_0 + I s_1\in C}$. By Corollary~\ref{PowersSimple}, we have $S_L^{-m}(s,x)\prodL S_L^{-n}(\overline{s},x)={ (s-x)^{-m}(\overline{s}-x)^{-n}}$ for any $x\in\C_I$. Lemma~\ref{DistLem} implies $\dist(s,\overline{U}\cap\C_{I}) = \dist(s,\overline{U}) > \varepsilon$. Since $\overline{U}\cap\C_I$ is symmetric with respect to the real axis, we also have $\dist(\overline{s},\overline{U}\cap\C_{I})  > \varepsilon$ and we deduce

\begin{align*}
\left\| S_L^{-m}(s,T)\prodL S_L^{-n}(\overline{s},T) \right\| &= \left\| \frac{1}{2\pi} \int_{\partial(U\cap\C_I)} S_L^{-1}(p,T)\,dp_I\, S_L^{-m}(s,p)\prodL S_L^{-n}(\overline{s},p) \right\|\\
&= \left\| \frac{1}{2\pi} \int_{\partial(U\cap\C_I)} S_L^{-1}(p,T)\,dp_I\, (s-p)^{-m}(\overline{s}-p)^{-n} \right\|\\
&\leq  \frac{1}{2\pi} \int_{\partial(U\cap\C_I)} \left\|S_L^{-1}(p,T)\right\|\,dp\, \left|(s-p)^{-m}(\overline{s}-p)^{-n}\right| \\
&\leq  \frac{1}{2\pi} \int_{\partial(U\cap\C_I)} \left\|S_L^{-1}(p,T)\right\|\,dp\, \frac{1}{\varepsilon^{m+n}}.
\end{align*}
 Hence, if we set
 \[K_T :=  \sup_{J\in\S}\frac{1}{2\pi}\int_{\partial(U\cap\C_J)} \left\|S_L^{-1}(p,T)\right\|\,dp\,,\]
 then
 \[ \left\| S_L^{-m}(s,T)\prodL S_L^{-n}(\overline{s},T)\right\| \leq \frac{K_T}{\varepsilon^{m+n}},\]
where the constant $K_T$ depends neither on the point $s\in C$ nor on the numbers $n,m\geq 0$.

\end{proof}

\begin{Tm}\label{SpectrumStable}
Let $T\in \boundOP^{0,1}(V_n)$ and let $N\in \boundOP^{0,1}(V_n)$ such that $T$ and $N$ commute and such that $\sigma_S(N)$ is contained in the open ball $B_{\varepsilon}(0)$.
If $\dist(s,\sigma_S(T)) > \varepsilon$, then $s\in\rho_S(T+N)$ and
\begin{equation*}
\left((T+N)^2 - 2s_0(T+N) + |s|^2\id\right)^{-1} = \sum_{n=0}^\infty\left(\sum_{k=0}^nS_L^{-(k+1)}(s,T)\prodL S_L^{-(n-k+1)}(\overline{s},T)\right)N^n,
\end{equation*}
where the series converges in the operator norm.
\end{Tm}
\begin{proof}
We first show the convergence of the series
\[\Sigma(s,T,N) := \sum_{n=0}^\infty\left(\sum_{k=0}^nS_L^{-(k+1)}(s,T)\prodL S_L^{-(n-k+1)}(\overline{s},T)\right)N^n.\]
 Since $\sigma_S(N)$ is compact, there exists $\theta\in(0,1)$ such that $\sigma_S(N) \subset B_{\theta\varepsilon}(0)\subset B_{\varepsilon}(0)$. Applying the $S$-functional calculus, we obtain
\begin{align*}
\| N^m\| &= \left\|\frac{1}{2\pi} \int_{\partial(B_{\theta\varepsilon}(0)\cap\C_I)}S_L^{-1}(s,N)\,ds_I\, s^m\right\|\\
&\leq \frac{1}{2\pi} \int_{\partial(B_{\theta\varepsilon}(0)\cap\C_I)}\left\|S_L^{-1}(s,N)\right\|\,ds\, |s|^m
\\
&
=\frac{1}{2\pi} \int_{\partial(B_{\theta\varepsilon}(0)\cap\C_I)}\left\|S_L^{-1}(s,N)\right\|\,ds\,  (\theta\varepsilon)^m
\end{align*}
for any $m\geq 0$. Hence,
\begin{equation}\label{K_NEstimate} \|N^m\| \leq K_N (\theta\varepsilon)^m\end{equation}
with
\[K_N:= \frac{1}{2\pi} \int_{\partial(B_{\theta\varepsilon}(0)\cap\C_I)}\left\|S_L^{-1}(s,N)\right\|\,ds.\]
From Proposition~\ref{KTEstimate}, we deduce
\begin{align*} & \sum_{n=0}^\infty\left\|\left(\sum_{k=0}^nS_L^{-(k+1)}(s,T)\prodL S_L^{-(n-k+1)}(\overline{s},T)\right)N^n \right\|\\
\leq& \sum_{n=0}^\infty \sum_{k=0}^n\left\|S_L^{-(k+1)}(s,T)\prodL S_L^{-(n-k+1)}(\overline{s},T) \right\| \left\|N^n \right\|\\
\leq& \sum_{n=0}^\infty(n+1) \frac{K_T}{\varepsilon^{n+2}} K_N(\theta\varepsilon)^n \leq \frac{K_TK_N}{\varepsilon^2} \sum_{n=0}^{\infty} (n+1)\theta^n.
\end{align*}
By the roots test, this last series converges because $0 < \theta < 1$. The comparison test yields the convergence of the original series $\Sigma(s,T,N)$ in the operator norm.

From Definition~\ref{ResSlicePowers} and the fact that $T$ and $N$ commute, we deduce
\begin{align*}
&\Sigma(s,T,N)((T+N)^2 - 2s_0(T+N) + |s|^2\id)\\
=&\left(\sum_{n=0}^{\infty}\left(\sum_{k=0}^{n}S_L^{-(k+1)}(s,T)\prodL S_L^{-(n-k+1)}(\overline{s},T)\right)N^n\right)((T^2-2TN + N^2 - 2s_0T-2s_0N + |s|^2\id)\\
=& \sum_{n=0}^\infty (T^2 - 2s_0T+|s|^2\id)^{-(n+2)}\left(\sum_{k=0}^n(\overline{s}\:\!\id - T)^{\prodL(k+1)}\prodL({s}\id-T)^{\prodL(n-k+1)}\right)N^n(T^2 - 2s_0 T + |s|^2\id)\\
&+\sum_{n=0}^\infty (T^2 - 2s_0T+|s|^2\id)^{-(n+2)}\left(\sum_{k=0}^n(\overline{s}\:\!\id - T)^{\prodL(k+1)}\prodL({s}\id-T)^{\prodL(n-k+1)}\right)N^n(2T- 2s_0\id) N\\
&+\sum_{n=0}^\infty (T^2 - 2s_0T+|s|^2\id)^{-(n+2)}\left(\sum_{k=0}^n(\overline{s}\:\!\id - T)^{\prodL(k+1)}\prodL({s}\id-T)^{\prodL(n-k+1)}\right)N^nN^2.
\end{align*}
Applying Corollary~\ref{PowersSimple} and the $S$-functional calculus, we see that any sum
\[\sum_{k=0}^n(\overline{s}\:\!\id - T)^{\prodL(k+1)}\prodL({s}\id-T)^{\prodL(n-k+1)}\] is a polynomial in $T$ with real coefficients and so it commutes with the operator ${T^2 -2s_0T +|s|^2\id}$. Remark~\ref{SeriesProdSimple} implies
\begin{align*}
&\Sigma(s,T,N)\,((T+N)^2 - 2s_0(T+N) + |s|^2\id)\\
=& \sum_{n=0}^\infty (T^2 - 2s_0T+|s|^2\id)^{-(n+1)}\left(\sum_{k=0}^n(\overline{s}\:\!\id - T)^{\prodL(k+1)}\prodL({s}\id-T)^{\prodL(n-k+1)}\right)N^n\\
&+\sum_{n=0}^\infty (T^2 - 2s_0T+|s|^2\id)^{-(n+2)}\left(\sum_{k=0}^n(\overline{s}\:\!\id - T)^{\prodL(k+1)}\prodL({s}\id-T)^{\prodL(n-k+1)}\prodL(2T- 2s_0\id)\right) N^{n+1}\\
&+\sum_{n=0}^\infty (T^2 - 2s_0T+|s|^2\id)^{-(n+2)}\left(\sum_{k=0}^n(\overline{s}\:\!\id - T)^{\prodL(k+1)}\prodL({s}\id-T)^{\prodL(n-k+1)}\right)N^{n+2}\displaybreak[1]\\
=& \sum_{n=0}^\infty (T^2 - 2s_0T+|s|^2\id)^{-(n+1)}\left(\sum_{k=0}^n(\overline{s}\:\!\id - T)^{\prodL(k+1)}\prodL({s}\id-T)^{\prodL(n-k+1)}\right)N^n\\
&+\sum_{n=1}^\infty (T^2 - 2s_0T+|s|^2\id)^{-(n+1)}\left(\sum_{k=0}^{n-1}(\overline{s}\:\!\id - T)^{\prodL(k+1)}\prodL({s}\id-T)^{\prodL(n-k)}\prodL(2T- 2s_0\id)\right) N^{n}\\
&+\sum_{n=2}^\infty (T^2 - 2s_0T+|s|^2\id)^{-n}\left(\sum_{k=0}^{n-2}(\overline{s}\:\!\id - T)^{\prodL(k+1)}\prodL({s}\id-T)^{\prodL(n-k-1)}\right)N^{n}.
\end{align*}
The identity
\begin{align*}&(T^2 - 2s_0T+|s|^2\id)^{-n}\left(\sum_{k=0}^{n-2}(\overline{s}\:\!\id - T)^{\prodL(k+1)}\prodL({s}\id-T)^{\prodL(n-k-1)}\right) \\
=& (T^2 - 2s_0T+|s|^2\id)^{-n}\left(\sum_{k=1}^{n-1}(\overline{s}\:\!\id - T)^{\prodL k}\prodL({s}\id-T)^{\prodL(n-k)}\right)\\
= &(T^2 - 2s_0T+|s|^2\id)^{-(n+1)}\left(\sum_{k=1}^{n-1}(\overline{s}\:\!\id - T)^{\prodL(k+1)}\prodL({s}\id-T)^{\prodL(n-k+1)}\right),\end{align*}
finally yields
\begin{align*}
&\Sigma(s,T,N)\,((T+N)^2 - 2s_0(T+N) + |s|^2\id)\\
=& \sum_{n=0}^\infty (T^2 - 2s_0T+|s|^2\id)^{-(n+1)}\Bigg(\sum_{k=0}^n(\overline{s}\:\!\id - T)^{\prodL(k+1)}\prodL({s}\id-T)^{\prodL(n-k+1)}\\
&\phantom{\sum_{n=0}^\infty (T^2 - 2s_0T+|s|^2\id)^{-(n+1)}}+\sum_{k=0}^{n-1}(\overline{s}\:\!\id - T)^{\prodL(k+1)}\prodL({s}\id-T)^{\prodL(n-k)}\prodL(2T- 2s_0\id)\\
&\phantom{\sum_{n=0}^\infty (T^2 - 2s_0T+|s|^2\id)^{-(n+1)}}+\sum_{k=1}^{n-1}(\overline{s}\:\!\id - T)^{\prodL(k+1)}\prodL({s}\id-T)^{\prodL(n-k+1)}\Bigg)N^{n},
\end{align*}
where the empty sum equals $0$. Hence,
\begin{gather*}\Sigma(s,T,N)\,((T+N)^2 - 2s_0(T+N) + |s|^2\id)
= \sum_{n=0}^\infty (T^2 - 2s_0T + |s|^2\id)^{-(n+1)}A_n N^n,\end{gather*}
with operator-valued coefficients
\begin{equation}\label{Coeffs}\begin{split}A_n =& \sum_{k=0}^n(\overline{s}\:\!\id - T)^{\prodL(k+1)}\prodL({s}\id-T)^{\prodL(n-k+1)}\\
&+\sum_{k=0}^{n-1}(\overline{s}\:\!\id - T)^{\prodL(k+1)}\prodL({s}\id-T)^{\prodL(n-k)}\prodL(2T- 2s_0\id)\\
&+\sum_{k=1}^{n-1}(\overline{s}\:\!\id - T)^{\prodL(k+1)}\prodL({s}\id-T)^{\prodL(n-k+1)}.
\end{split}
\end{equation}
We have
\[A_0 =  (\overline{s}\:\!\id - T)\prodL({s}\id-T)^{\prodL 1}=T^2 - 2 s_0 + |s|^2\id\]
and
\begin{align*} A_1 &= (\overline{s}\:\!\id - T)\prodL({s}\id-T)^{\prodL 2}+(\overline{s}\:\!\id - T)^{\prodL 2}\prodL({s}\id-T)+(\overline{s}\:\!\id - T)^{\prodL 1}\prodL({s}\:\id-T)^{\prodL 1}\prodL(2T- 2s_0\id)\\
&= (\overline{s}\:\!\id - T)\prodL({s}\id-T)\prodL({s}\id-T -\overline{s}\:\!\id - T)+(\overline{s}\:\!\id - T)\prodL({s}\id-T)\prodL(2T- 2s_0\id) = 0,
\end{align*}
where we used  that
\[{(\overline{s}\:\!\id - T)\prodL(s\id-T)=}  {T^2 - 2s_0 T + |s|^2\id =} {(s\id-T)\prodL(\overline{s}\:\!\id - T)}.\] By adding the terms for  $k=0,\ldots, n-1$ of the fist sum in \eqref{Coeffs} to the second sum and the last term with $k = n$ to the third sum, we obtain for $n\geq 2$

 \begin{align*}
A_n =& \sum_{k=0}^{n-1}(\overline{s}\:\!\id - T)^{\prodL(k+1)}\prodL({s}\id-T)^{\prodL(n-k)}\prodL({s}\id - T + 2T- 2s_0\id)\\
&+\sum_{k=1}^{n}(\overline{s}\:\!\id - T)^{\prodL(k+1)}\prodL({s}\id-T)^{\prodL(n-k+1)}\\
=& \sum_{k=0}^{n-1}(\overline{s}\:\!\id - T)^{\prodL(k+1)}\prodL({s}\id-T)^{\prodL(n-k)}\prodL(  T- \overline{s}\:\!\id) \\
&+\sum_{k=0}^{n-1}(\overline{s}\:\!\id - T)^{\prodL(k+1)}\prodL(s\id-T)^{\prodL(n-k)}\prodL(\overline{s}\:\!\id - T) = 0.
\end{align*}
Therefore,
\begin{gather*}\Sigma(s,T,N)\,((T+N)^2 - 2s_0(T+N) + |s|^2\id) \\
= (T^2 - 2s_0T + |s|^2\id)^{-1} (T^2 - 2s_0T + |s|^2\id)^{-1} = \id.\end{gather*}
The coefficients of the series $\Sigma(s,T,N)$ are
\begin{multline*}\sum_{k=0}^nS_L^{-(k+1)}(s,T)\prodL S_L^{-(n-k+1)}(\overline{s},T) \\
= (T^2 - 2s_0T + |s|^2\id)^{-(n+2)}\sum_{k=0}^n(\overline{s}\:\!\id-T)^{\prodL(k+1)}\prodL(s\id - T)^{\prodL(n-k+1)}.\end{multline*}
Corollary~\ref{PowersSimple} and the $S$-functional calculus imply that each sum
\[\sum_{k=0}^n(\overline{s}\:\!\id-T)^{\prodL(k+1)}\prodL(s\id - T)^{\prodL(n-k+1)}\]
is a polynomial in $T$ with real coefficients. Therefore, it commutes with $T$ and $N$. Similarly, also the operator $(T^2 - 2s_0T + |s|^2\id)^{-1}$ commutes with $T$ and $N$. Consequently, each coefficient of $\Sigma(s,T,N)$, and therefore the entire series $\Sigma(s,T,N)$, commutes with $T$ and $N$, and in turn, also with the operator ${(T+N)^2 + 2s_0 (T+N) + |s|^2\id}$. Hence, we also have
\[\left( (T+N)^2 -2s_0(T+N)+|s|^2\id\right) \Sigma(s,T,N)= \id,\]
and the operator $(T+N)^2 -2s_0(T+N)+|s|^2\id$ is invertible, which implies $s\in\rho_S(T+N)$.

\end{proof}

\begin{samepage}
\begin{Tm} \label{ResolventTaylor}
Let $T\in\boundOP^{0,1}(V_n)$ and let $N\in\boundOP^{0,1}(V_n)$ be such that $\sigma_S(N)\subset B_{\varepsilon}(0)$. For any $s\in\rho_S(T)$ with $\dist(s,\sigma_S(T)) > \varepsilon$, the identities
\[S_L^{-1}(s,T+N) = \sum_{n=0}^\infty N^n\, S_L^{-(n+1)}(s,T) \qquad\text{and}\qquad S_R^{-1}(s,T+N) = \sum_{n=0}^{\infty}S_R^{-(n+1)}(s,T)\, N^n\]
hold true. Moreover, these series converge uniformly on any set $C$ with ${\dist(C,\sigma_S(T))>\varepsilon}$.
\end{Tm}
\end{samepage}
\begin{proof}
In \eqref{K_NEstimate}, we showed the existence of two constants $K_N\geq 0$ and $\theta\in (0,1)$ such that
$\|N\|^m \leq K_N(\theta\varepsilon)^m$
for any $m\in\N_0$. Moreover, for any $C\subset \R^{n+1}$ with ${\dist(C,\sigma_S(T)) > \varepsilon}$, Proposition~\ref{KTEstimate} implies the existence of a constant $K_T$ such that $\|S_L^{-m}(s,T) \|\leq K_T/\varepsilon^m $ for any $s\in C$ and any $m\in\N_0$. Therefore, the estimate
\begin{align*} \sum_{n=n_0}^\infty \left\| N^n\,S_L^{-(n+1)}(s,T)\right\| &\leq \sum_{n=n_0}^\infty \left\| N^n\right\|\left\|S_L^{-(n+1)}(s,T)\right\| \\
&\leq  \sum_{n=n_0}^{\infty}K_N(\theta\varepsilon)^n\frac{K_T}{\varepsilon^{n+1}}
\\
&
= \frac{K_TK_N}{\varepsilon}\sum_{n=n_0}^\infty \theta^n \overset{n_0\to\infty}{\longrightarrow}0
\end{align*}
 holds true for any $s\in C$ and implies the uniform convergence of the series on $C$.

Let $s\in\rho_S(T)$ with $\dist(s,\sigma_S(T)) > \varepsilon$. We have
\begin{align*}
&\left((T+N)^2 - 2s_0(T+N)+|s|^2\id\right) \sum_{n=0}^{\infty}N^n  \,S_L^{-(n+1)}(s,T)
 \\
=&\left(T^2 - 2s_0T+|s|^2\id\right) \sum_{n=0}^{\infty}N^n  (T^2 - 2s_0T + |s|^2\id)^{-(n+1)}(\overline{s}\:\!\id-T)^{\prodL(n+1)}\\
&+\left( 2T-2s_0\id\right)N \sum_{n=0}^{\infty}N^n (T^2 - 2s_0T + |s|^2\id)^{-(n+1)}(\overline{s}\:\!\id-T)^{\prodL(n+1)}\\
&+N^2 \sum_{n=0}^{\infty}N^n (T^2 - 2s_0T + |s|^2\id)^{-(n+1)}(\overline{s}\:\!\id-T)^{\prodL(n+1)}\\
\displaybreak[1]
=&\sum_{n=0}^{\infty}N^n  (T^2 - 2s_0T + |s|^2\id)^{-n}(\overline{s}\:\!\id-T)^{\prodL(n+1)}\\
&+\sum_{n=0}^{\infty}N^{n+1}\left( 2T-2s_0\id\right)(T^2 - 2s_0T + |s|^2\id)^{-(n+1)}(\overline{s}\:\!\id-T)^{\prodL(n+1)}\\
&+\sum_{n=0}^{\infty}N^{n+2} (T^2 - 2s_0T + |s|^2\id)^{-(n+1)}(\overline{s}\:\!\id-T)^{\prodL(n+1)}.
\end{align*}
Shifting the indices yields
\begin{align*}
&\left((T+N)^2 - 2s_0(T+N)+|s|^2\id\right) \sum_{n=0}^{\infty}N^n  \,S_L^{-(n+1)}(s,T) \\
=&\sum_{n=0}^{\infty}N^n  (T^2 - 2s_0T + |s|^2\id)^{-n}(\overline{s}\:\!\id-T)^{\prodL(n+1)}\\
&+\sum_{n=1}^{\infty}N^{n}\left(2T- 2s_0\id\right) (T^2 - 2s_0T + |s|^2\id)^{-n}(\overline{s}\:\!\id-T)^{\prodL n}\\
&+\sum_{n=2}^{\infty}N^{n} (T^2 - 2s_0T + |s|^2\id)^{-(n-1)}(\overline{s}\:\!\id-T)^{\prodL(n-1)}\displaybreak[1]\\
=& \overline{s}\:\!\id-T + N(T^2 -s_0 T + |s|^2\id)^{-1}(\overline{s}\:\!\id -T)^{\prodL 2}
\\
&
+ N(T^2 - 2s_0T +|s|^2\id)^{-1}(2T-2s_0\id)(\overline{s}\:\!\id-T)+\\
&+ \sum_{n=2}^{\infty}N^n (T^2 - 2s_0T + |s|^2\id)^{-n}\Big[(\overline{s}\:\!\id-T)^{\prodL(n+1)}+\left(2T-2s_0\id \right) (\overline{s}\:\!\id-T)^{\prodL n}\\
&\qquad\qquad+ (T^2 - 2s_0T + |s|^2\id)(\overline{s}\:\!\id-T)^{\prodL(n-1)}\Big].
\end{align*}
The last series equals $0$ because Remark~\ref{SeriesProdSimple} and the identity
$$
(T^2 - 2s_0 T + |s|^2\id) = (s\id- T)\prodL (\overline{s}\:\!\id-T)
$$
imply
\begin{align*}
&(\overline{s}\:\!\id-T)^{\prodL(n+1)}+\left(2 T-2s_0\id\right) (\overline{s}\:\!\id-T)^{\prodL n}+ (T^2 - 2s_0T + |s|^2\id)(\overline{s}\:\!\id-T)^{\prodL(n-1)}\\
=&(\overline{s}\:\!\id-T)^{\prodL(n+1)}+\left(2 T-2s_0\id\right)\prodL (\overline{s}\:\!\id-T)^{\prodL n}+ (s\id - T)\prodL(\overline{s}\:\!\id-T)^{\prodL n}\\
=& (\overline{s}\:\!\id - T + 2T - 2s_0\id + s\id- T)\prodL(\overline{s}\:\!\id - T)^{\prodL(n-1)} = 0.
\end{align*}
Hence, we finally obtain
 \begin{align*}&\left((T+N)^2 - 2s_0(T+N)+|s|^2\id\right) \sum_{n=0}^{\infty}N^n  S_L^{-(n+1)}(s,T)
  \\
=& \overline{s}\:\!\id-T + N(T^2 -2s_0 T + |s|^2\id)^{-1}(\overline{s}^2\id -  2T\overline{s} + T^2) \\
&+ N(T^2 - 2s_0T +|s|^2\id)^{-1}(2T\overline{s} - 2s_0\overline{s}\:\!\id - 2T^2 + 2s_0T) \\
=&\overline{s} \id -T + N(T^2 - 2s_0 T + |s|^2\id)^{-1}(-T^2 + 2s_0T - |s|^2\id) = \overline{s}\:\!\id-T-N.
  \end{align*}
  Since $(T+N)^2-2s_0(T+N) + |s|^2\id$ is invertible by Theorem~\ref{SpectrumStable}, this is equivalent to
    \[ \sum_{n=0}^{\infty}N^n \, S_L^{-(n+1)}(s,T) = ((T+N)^2 - 2s_0(T+N) + |s|^2\id)^{-1}(\overline{s}\:\!\id-T-N) = S_L^{-1}(s,T+N).\]

\end{proof}

In order to show the Taylor formula in the operator variable for the $S$-functional calculus, we calculate the slice derivatives of the $S$-resolvent operators with respect to the scalar variable.

\begin{Pn}[Derivatives of the $S$-resolvent operators]\label{ResDeriv}
Let  $T\in \mathcal{B}^{\small 0,1}(V_n)$ and let ${s\in \rho_S(T)}$.
Then
\begin{equation}\label{quatSresolrddlft}
\partial_s^{m}S_L^{-1}(s,T)= (-1)^m m!\,S_L^{-(m+1)}(s,T)
\end{equation}
and
\begin{equation}\label{quatSresorig}
\partial_s^{m}S_R^{-1}(s,T)=(-1)^mm!\,S_R^{-(m+1)}(s,T),
\end{equation}
for any $m\geq 0$.
\end{Pn}
\begin{proof}
Recall from Remark~\ref{SliceDerivS0} that, for slice hyperholomorphic functions, the slice derivative defined in Definition~\ref{SliceDeriv} coincides with
the partial derivative with respect to the real part $s_0$. We show only \eqref{quatSresolrddlft}, since \eqref{quatSresorig} follows by analogous computations.

We prove the statement by induction. For $m=0$, the identity \eqref{quatSresolrddlft} is obvious. We assume that
$\partial_{s_0}^{m-1}S_L^{-1}(s,T)=(-1)^{m-1}(m-1)!\,S_L^{-m}(s,T)$ and we compute $\partial_{s_0}^{m}S_L^{-1}(s,T)$. Let $I_s\in\S$ such that $s\in\C_{I_s}$. If $x\in\C_{I_s}$, then  $S_L^{-n}(s,x)  = (s-x)^{-n}$ by Corollary~\ref{PowersSimple}.
Applying the $S$-functional calculus yields
\[\begin{split}
\frac{\partial}{\partial s_0} S_L^{-m}(s,T) &= \frac{\partial}{\partial s_0} \frac{1}{2\pi}\int_{\partial(U\cap\C_{I_s})} S_L^{-1}(p,T)\,dp_I\,(s-p)^{-m}
\\
&= \frac{1}{2\pi}\int_{\partial(U\cap\C_{I_s})} S_L^{-1}(p,T)\,dp_I\,\frac{\partial}{\partial s_0} (s-p)^{-m}
\\
&= -m \frac{1}{2\pi}\int_{\partial(U\cap\C_{I_s})} S_L^{-1}(p,T)\,dp_I\,(s-p)^{-(m+1)}
\\
&
= -m\,  S_L^{-(m+1)}(s,T),
\end{split}\]
and in turn,
\[\partial_s^m S_L^{-1}(s,T) = \frac{\partial^m}{\partial s_0^m} S_L^{-1}(s,T) = \frac{\partial}{\partial s_0} (-1)^{m-1}(m-1)!\,S_L^{-m}(s,T) = (-1)^mm!\, S_L^{-(m+1)}(s,T).\]

\end{proof}

\begin{Tm}[The Taylor formulas]\label{THTAYLOR}
Let $T\in\boundOP^{0,1}(V_n)$  and $N\in\boundOP^{0,1}(V_n)$ with ${\sigma_S(N)\subset B_{\varepsilon}(0)}$. If $f\in \mathcal{SM}^L(U)$ for some $T$-admissible domain $U$ that contains $\sigma_S(T)$  and any point $s\in\rr^{n+1}$ with $\dist( s,\sigma_S(T))\leq \varepsilon$, then $f\in \mathcal{SM}^L_{\sigma_S(T+N)}$ and
\[
f(T+N) =  \sum_{n=0}^{\infty} N^n\frac{1}{n!} \left(\partial_{s}^nf\right)(T).
 \]
Similarly, if $f\in \mathcal{SM}^R(U)$, for some $T$-admissible domain $U$ that contains $\sigma_S(T)$  and any point $s\in\H$ with $\dist( s,\sigma_S(T))\leq \varepsilon$, then $f\in\SM^R_{\sigma_S(T+N)}$ and
\[
f(T+N) = \sum_{n=0}^\infty \frac{1}{n!}\left(\partial_{s}^nf\right)(T)N^n.
\]
\end{Tm}
\begin{proof}
We prove just the first Taylor formula, the second one is obtained with similar computations.
It follows from Theorem~\ref{SpectrumStable} that $\sigma_S(T+N)\subset U$, and hence, $f\in\SM^L_{\sigma_S(T+N)}$. Applying the $S$-functional calculus and Theorem~\ref{ResolventTaylor}, we obtain
\[
\begin{split} f(T+N) &= \frac{1}{2\pi}\int_{\partial (U\cap\mathbb C_I)}S_L^{-1}(s,T+N) \,ds_I\, f(s)
\\
&
= \frac{1}{2\pi}\int_{\partial (U\cap\mathbb C_I)} \sum_{n= 0}^{\infty}N^nS_L^{-(n+1)}(s,T)\,ds_I\, f(s)
\\
&
=  \sum_{n= 0}^{\infty}N^n\frac{1}{2\pi}\int_{\partial (U\cap\mathbb C_I)}S_L^{-(n+1)}(s,T)\,ds_I\, f(s).
\end{split}
\]
By Proposition~\ref{ResDeriv}, we have
$$
S_L^{-(n+1)}(s,T)  = (-1)^n \frac{1}{n!}\partial_{s}^n S_L^{-1}(s,T).
$$
 Hence,
\[f(T+N) = \sum_{n=0}^{\infty} N^n \frac{(-1)^n}{n!} \frac{1}{2\pi}\int_{\partial (U\cap\mathbb C_I)} \partial_{s}^nS_L^{-1}(s,T)   \,ds_I\, f(s) . \]
By integrating the $n$-th term in the sum $n$-times by parts, we finally obtain
\[f(T+N) = \sum_{n= 0}^{\infty}N^n \frac{1}{n!}\ \frac{1}{2\pi}\int_{\partial (U\cap\mathbb C_I)} S_L^{-1}(s,T) \,ds_I\, (\partial_{s}^nf)(s) = \sum_{n=0}^{\infty} N^n\frac{1}{n!} (\partial_{s}^nf)(T). \]

\end{proof}

\begin{Rk}
In Proposition \ref{KTEstimate}
it can be  C=$\R$. Then we cannot use a slice domain,
but this limitation is removed when we  use the definition of slice monogenic functions as in Remark~\ref{RKTWODEF} which requires
just axially symmetric domains.
\end{Rk}

\end{document}